\documentclass[12pt,a4paper]{amsart}
\usepackage{amsmath,amssymb,amscd}

\makeatletter
\def\diagram{\leftwidth=\z@ \rightwidth=\z@ \topheight=\z@
\botheight=\z@ \setbox\@picbox\hbox\bgroup}
\def\enddiagram{\egroup\wd\@picbox\rightwidth\unitlength
\ht\@picbox\topheight\unitlength \dp\@picbox\botheight\unitlength
\hskip\leftwidth\unitlength\box\@picbox}
\def\bfig{\begin{diagram}}
\def\efig{\end{diagram}}
\newcount\wideness \newcount\leftwidth \newcount\rightwidth
\newcount\highness \newcount\topheight \newcount\botheight
\def\ratchet#1#2{\ifnum#1<#2 \global #1=#2 \fi}
\def\putbox(#1,#2)#3{%
\horsize{\wideness}{#3} \divide\wideness by 2 {\advance\wideness
by #1 \ratchet{\rightwidth}{\wideness}} {\advance\wideness by -#1
\ratchet{\leftwidth}{\wideness}} \vertsize{\highness}{#3}
\divide\highness by 2 {\advance\highness by #2
\ratchet{\topheight}{\highness}} {\advance\highness by -#2
\ratchet{\botheight}{\highness}} \put(#1,#2){\makebox(0,0){$#3$}}}
\def\putlbox(#1,#2)#3{%
\horsize{\wideness}{#3} {\advance\wideness by #1
\ratchet{\rightwidth}{\wideness}} {\ratchet{\leftwidth}{-#1}}
\vertsize{\highness}{#3} \divide\highness by 2 {\advance\highness
by #2 \ratchet{\topheight}{\highness}} {\advance\highness by -#2
\ratchet{\botheight}{\highness}}
\put(#1,#2){\makebox(0,0)[l]{$#3$}}}
\def\putrbox(#1,#2)#3{%
\horsize{\wideness}{#3} {\ratchet{\rightwidth}{#1}}
{\advance\wideness by -#1 \ratchet{\leftwidth}{\wideness}}
\vertsize{\highness}{#3} \divide\highness by 2 {\advance\highness
by #2 \ratchet{\topheight}{\highness}} {\advance\highness by -#2
\ratchet{\botheight}{\highness}}
\put(#1,#2){\makebox(0,0)[r]{$#3$}}}

\def\adjust[#1]{} 
\newcount \coefa
\newcount \coefb
\newcount \coefc
\newcount\tempcounta
\newcount\tempcountb
\newcount\tempcountc
\newcount\tempcountd
\newcount\xext
\newcount\yext
\newcount\xoff
\newcount\yoff
\newcount\gap%
\newcount\arrowtypea
\newcount\arrowtypeb
\newcount\arrowtypec
\newcount\arrowtyped
\newcount\arrowtypee
\newcount\height
\newcount\width
\newcount\xpos
\newcount\ypos
\newcount\run
\newcount\rise
\newcount\arrowlength
\newcount\halflength
\newcount\arrowtype
\newdimen\tempdimen
\newdimen\xlen
\newdimen\ylen
\newsavebox{\tempboxa}%
\newsavebox{\tempboxb}%
\newsavebox{\tempboxc}%
\newdimen\w@dth
\def\setw@dth#1#2{\setbox\z@\hbox{$#1$}\w@dth=\wd\z@
\setbox\@ne\hbox{$#2$}\ifnum\w@dth<\wd\@ne \w@dth=\wd\@ne \fi
\advance\w@dth by 1.2em}
\def\t@^#1_#2{\def\n@one{#1}\def\n@two{#2}\mathrel{\setw@dth{#1}{#2}
\mathop{\hbox to \w@dth{\rightarrowfill}}\limits
\ifx\n@one\empty\else ^{\box\z@}\fi \ifx\n@two\empty\else
_{\box\@ne}\fi}}
\def\t@@^#1{\@ifnextchar_ {\t@^{#1}}{\t@^{#1}_{}}}
\def\to{\@ifnextchar^ {\t@@}{\t@@^{}}}
\def\t@left^#1_#2{\def\n@one{#1}\def\n@two{#2}\mathrel{\setw@dth{#1}{#2}
\mathop{\hbox to \w@dth{\leftarrowfill}}\limits
\ifx\n@one\empty\else ^{\box\z@}\fi \ifx\n@two\empty\else
_{\box\@ne}\fi}}
\def\t@@left^#1{\@ifnextchar_ {\t@left^{#1}}{\t@left^{#1}_{}}}
\def\toleft{\@ifnextchar^ {\t@@left}{\t@@left^{}}}
\def\two@^#1_#2{\def\n@one{#1}\def\n@two{#2}\mathrel{\setw@dth{#1}{#2}
\mathop{\vcenter{\hbox to \w@dth{\rightarrowfill}\kern-1.7ex
                 \hbox to \w@dth{\rightarrowfill}}%
       }\limits
\ifx\n@one\empty\else ^{\box\z@}\fi \ifx\n@two\empty\else
_{\box\@ne}\fi}}
\def\tw@@^#1{\@ifnextchar_ {\two@^{#1}}{\two@^{#1}_{}}}
\def\two{\@ifnextchar^ {\tw@@}{\tw@@^{}}}
\def\tofr@^#1_#2{\def\n@one{#1}\def\n@two{#2}\mathrel{\setw@dth{#1}{#2}
\mathop{\vcenter{\hbox to \w@dth{\rightarrowfill}\kern-1.7ex
                 \hbox to \w@dth{\leftarrowfill}}%
       }\limits
\ifx\n@one\empty\else ^{\box\z@}\fi \ifx\n@two\empty\else
_{\box\@ne}\fi}}
\def\t@fr@^#1{\@ifnextchar_ {\tofr@^{#1}}{\tofr@^{#1}_{}}}
\def\tofro{\@ifnextchar^ {\t@fr@}{\t@fr@^{}}}

\def\mon{\mathop{\m@th\hbox to
      14.6\P@{\lasyb\char'51\hskip-2.1\P@$\arrext$\hss
$\mathord\rightarrow$}}\limits} 
\def\leftmono{\mathrel{\m@th\hbox to
14.6\P@{$\mathord\leftarrow$\hss$\arrext$\hskip-2.1\P@\lasyb\char'50%
}}\limits} 
\mathchardef\arrext="0200       

\def\settypes(#1,#2,#3){\arrowtypea#1 \arrowtypeb#2 \arrowtypec#3}
\def\settoheight#1#2{\setbox\@tempboxa\hbox{#2}#1\ht\@tempboxa\relax}%
\def\settodepth#1#2{\setbox\@tempboxa\hbox{#2}#1\dp\@tempboxa\relax}%
\def\settokens[#1`#2`#3`#4]{%
     \def\tokena{#1}\def\tokenb{#2}\def\tokenc{#3}\def\tokend{#4}}
\def\setsqparms[#1`#2`#3`#4;#5`#6]{%
\arrowtypea #1 \arrowtypeb #2 \arrowtypec #3 \arrowtyped #4
\width #5 \height #6 }
\def\setpos(#1,#2){\xpos=#1 \ypos#2}

\def\settriparms[#1`#2`#3;#4]{\settripairparms[#1`#2`#3`1`1;#4]}%
\def\settripairparms[#1`#2`#3`#4`#5;#6]{%
\arrowtypea #1 \arrowtypeb #2 \arrowtypec #3 \arrowtyped #4
\arrowtypee #5 \width #6 \height #6 }
\def\resetparms{\settripairparms[1`1`1`1`1;500]\width 500}
\resetparms
\def\mvector(#1,#2)#3{
\put(0,0){\vector(#1,#2){#3}}%
\put(0,0){\vector(#1,#2){26}}%
}
\def\evector(#1,#2)#3{{
\arrowlength #3
\put(0,0){\vector(#1,#2){\arrowlength}}%
\advance \arrowlength by-30
\put(0,0){\vector(#1,#2){\arrowlength}}%
}}
\def\horsize#1#2{%
\settowidth{\tempdimen}{$#2$}%
#1=\tempdimen \divide #1 by\unitlength }
\def\vertsize#1#2{%
\settoheight{\tempdimen}{$#2$}%
#1=\tempdimen
\settodepth{\tempdimen}{$#2$}%
\advance #1 by\tempdimen \divide #1 by\unitlength }
\def\putvector(#1,#2)(#3,#4)#5#6{{%
\ifnum3<\arrowtype \putdashvector(#1,#2)(#3,#4)#5\arrowtype \else
\ifnum\arrowtype<-3 \putdashvector(#1,#2)(#3,#4)#5\arrowtype \else
\xpos=#1 \ypos=#2 \run=#3 \rise=#4 \arrowlength=#5 \ifnum
\arrowtype<0
    \ifnum \run=0
        \advance \ypos by-\arrowlength
    \else
        \tempcounta \arrowlength
        \multiply \tempcounta by\rise
        \divide \tempcounta by\run
        \ifnum\run>0
            \advance \xpos by\arrowlength
            \advance \ypos by\tempcounta
        \else
            \advance \xpos by-\arrowlength
            \advance \ypos by-\tempcounta
        \fi
    \fi
    \multiply \arrowtype by-1
    \multiply \rise by-1
    \multiply \run by-1
\fi \ifcase \arrowtype
\or \put(\xpos,\ypos){\vector(\run,\rise){\arrowlength}}%
\or \put(\xpos,\ypos){\mvector(\run,\rise)\arrowlength}%
\or \put(\xpos,\ypos){\evector(\run,\rise){\arrowlength}}%
\fi\fi\fi }}
\def\putsplitvector(#1,#2)#3#4{
\xpos #1 \ypos #2 \arrowtype #4 \halflength #3 \arrowlength #3
\gap 140 \advance \halflength by-\gap \divide \halflength by2
\ifnum\arrowtype>0
   \ifcase \arrowtype
   \or \put(\xpos,\ypos){\line(0,-1){\halflength}}%
       \advance\ypos by-\halflength
       \advance\ypos by-\gap
       \put(\xpos,\ypos){\vector(0,-1){\halflength}}%
   \or \put(\xpos,\ypos){\line(0,-1)\halflength}%
       \put(\xpos,\ypos){\vector(0,-1)3}%
       \advance\ypos by-\halflength
       \advance\ypos by-\gap
       \put(\xpos,\ypos){\vector(0,-1){\halflength}}%
   \or \put(\xpos,\ypos){\line(0,-1)\halflength}%
       \advance\ypos by-\halflength
       \advance\ypos by-\gap
       \put(\xpos,\ypos){\evector(0,-1){\halflength}}%
   \fi
\else \arrowtype=-\arrowtype
   \ifcase\arrowtype
   \or \advance \ypos by-\arrowlength
       \put(\xpos,\ypos){\line(0,1){\halflength}}%
       \advance\ypos by\halflength
       \advance\ypos by\gap
       \put(\xpos,\ypos){\vector(0,1){\halflength}}%
   \or \advance \ypos by-\arrowlength
       \put(\xpos,\ypos){\line(0,1)\halflength}%
       \put(\xpos,\ypos){\vector(0,1)3}%
       \advance\ypos by\halflength
       \advance\ypos by\gap
       \put(\xpos,\ypos){\vector(0,1){\halflength}}%
   \or \advance \ypos by-\arrowlength
       \put(\xpos,\ypos){\line(0,1)\halflength}%
       \advance\ypos by\halflength
       \advance\ypos by\gap
       \put(\xpos,\ypos){\evector(0,1){\halflength}}%
   \fi
\fi }
\def\putmorphism(#1)(#2,#3)[#4`#5`#6]#7#8#9{{%
\run #2 \rise #3 \ifnum\rise=0
  \puthmorphism(#1)[#4`#5`#6]{#7}{#8}#9%
\else\ifnum\run=0
  \putvmorphism(#1)[#4`#5`#6]{#7}{#8}#9%
\else
\setpos(#1)%
\arrowlength #7 \arrowtype #8 \ifnum\run=0 \else\ifnum\rise=0
\else \ifnum\run>0
    \coefa=1
\else
   \coefa=-1
\fi \ifnum\arrowtype>0
   \coefb=0
   \coefc=-1
\else
   \coefb=\coefa
   \coefc=1
   \arrowtype=-\arrowtype
\fi \width=2 \multiply \width by\run \divide \width by\rise
\ifnum \width<0  \width=-\width\fi \advance\width by60 \if l#9
\width=-\width\fi
\putbox(\xpos,\ypos){#4}
{\multiply \coefa by\arrowlength
\advance\xpos by\coefa \multiply \coefa by\rise \divide \coefa
by\run \advance \ypos by\coefa
\putbox(\xpos,\ypos){#5} }%
{\multiply \coefa by\arrowlength
\divide \coefa by2 \advance \xpos by\coefa \advance \xpos by\width
\multiply \coefa by\rise \divide \coefa by\run \advance \ypos
by\coefa
\if l#9%
   \putrbox(\xpos,\ypos){#6}%
\else\if r#9%
   \putlbox(\xpos,\ypos){#6}%
\fi\fi }%
{\multiply \rise by-\coefc
\multiply \run by-\coefc \multiply \coefb by\arrowlength \advance
\xpos by\coefb \multiply \coefb by\rise \divide \coefb by\run
\advance \ypos by\coefb \multiply \coefc by70 \advance \ypos
by\coefc \multiply \coefc by\run \divide \coefc by\rise \advance
\xpos by\coefc \multiply \coefa by140 \multiply \coefa by\run
\divide \coefa by\rise \advance \arrowlength by\coefa
\ifcase\arrowtype
\or \put(\xpos,\ypos){\vector(\run,\rise){\arrowlength}}%
\or \put(\xpos,\ypos){\mvector(\run,\rise){\arrowlength}}%
\or \put(\xpos,\ypos){\evector(\run,\rise){\arrowlength}}%
\fi}\fi\fi\fi\fi}}

\newcount\numbdashes \newcount\lengthdash \newcount\increment
\def\howmanydashes{
\numbdashes=\arrowlength \lengthdash=40 \divide\numbdashes by
\lengthdash \lengthdash=\arrowlength \divide\lengthdash by
\numbdashes
\increment=\lengthdash \multiply\lengthdash by 3
\divide\lengthdash by 5 }
\def\putdashvector(#1)(#2,#3)#4#5{%
\ifnum#3=0 \putdashhvector(#1){#4}#5 \else \ifnum#2=0
\putdashvvector(#1){#4}#5\fi\fi}
\def\putdashhvector(#1,#2)#3#4{{%
\arrowlength=#3 \howmanydashes
\multiput(#1,#2)(\increment,0){\numbdashes}%
{\vrule height .4pt width \lengthdash\unitlength} \arrowtype=#4
\xpos=#1 \ifnum\arrowtype<0 \advance\arrowtype by 7 \fi
\ifcase\arrowtype \or \advance\xpos by 10
    \put(\xpos,#2){\vector(-1,0){\lengthdash}}
    \advance\xpos by 40
    \put(\xpos,#2){\vector(-1,0){\lengthdash}}
\or \advance \xpos by 10
    \put(\xpos,#2){\vector(-1,0){\lengthdash}}
    \advance\xpos by  \arrowlength
    \advance\xpos by  -50
    \put(\xpos,#2){\vector(-1,0){\lengthdash}}
\or \advance\xpos by 10
    \put(\xpos,#2){\vector(-1,0){\lengthdash}}
\or \advance\xpos by \arrowlength
    \advance\xpos by -\lengthdash
    \put(\xpos,#2){\vector(1,0){\lengthdash}}
\or {\advance\xpos by 10
    \put(\xpos,#2){\vector(1,0){\lengthdash}}}
    \advance\xpos by \arrowlength
    \advance\xpos by -\lengthdash
    \put(\xpos,#2){\vector(1,0){\lengthdash}}
\or \advance\xpos by \arrowlength
    \advance\xpos by -\lengthdash
    \put(\xpos,#2){\vector(1,0){\lengthdash}}
    \advance\xpos by -40
    \put(\xpos,#2){\vector(1,0){\lengthdash}}
   \fi
}}
\def\putdashvvector(#1,#2)#3#4{{%
\arrowlength=#3 \howmanydashes \ypos=#2 \advance\ypos by
-\arrowlength
\multiput(#1,#2)(0,\increment){\numbdashes}%
    {\vrule width .4pt height \lengthdash\unitlength}
\arrowtype=#4 \ypos=#2 \ifnum\arrowtype<0 \advance\arrowtype by 7
\fi \ifcase\arrowtype \or \advance\ypos by \arrowlength
\advance\ypos by -40
    \put(#1,\ypos){\vector(0,1){\lengthdash}}
    \advance\ypos by -40
    \put(#1,\ypos){\vector(0,1){\lengthdash}}
\or \advance\ypos by 10
    \put(#1,\ypos){\vector(0,1){\lengthdash}}
    \advance\ypos by \arrowlength \advance\ypos by -40
    \put(#1,\ypos){\vector(0,1){\lengthdash}}
\or \advance\ypos by \arrowlength \advance\ypos by -40
    \put(#1,\ypos){\vector(0,1){\lengthdash}}
\or \advance\ypos by 10
    \put(#1,\ypos){\vector(0,-1){\lengthdash}}
\or \advance\ypos by 10
    \put(#1,\ypos){\vector(0,-1){\lengthdash}}
    \advance\ypos by \arrowlength \advance\ypos by -40
    \put(#1,\ypos){\vector(0,-1){\lengthdash}}
\or \advance\ypos by 10
    \put(#1,\ypos){\vector(0,-1){\lengthdash}}
    \advance\ypos by 40
    \put(#1,\ypos){\vector(0,-1){\lengthdash}}
\fi }}
\def\puthmorphism(#1,#2)[#3`#4`#5]#6#7#8{{%
\xpos #1 \ypos #2 \width #6 \arrowlength #6 \arrowtype=#7
\putbox(\xpos,\ypos){#3\vphantom{#4}}%
{\advance \xpos by\arrowlength
\putbox(\xpos,\ypos){\vphantom{#3}#4}}%
\horsize{\tempcounta}{#3}%
\horsize{\tempcountb}{#4}%
\divide \tempcounta by2 \divide \tempcountb by2 \advance
\tempcounta by30 \advance \tempcountb by30 \advance \xpos
by\tempcounta \advance \arrowlength by-\tempcounta \advance
\arrowlength by-\tempcountb
\putvector(\xpos,\ypos)(1,0)\arrowlength\arrowtype \divide
\arrowlength by2 \advance \xpos by\arrowlength
\vertsize{\tempcounta}{#5}%
\divide\tempcounta by2 \advance \tempcounta by20
\if a#8 %
   \advance \ypos by\tempcounta
   \putbox(\xpos,\ypos){#5}%
\else
   \advance \ypos by-\tempcounta
   \putbox(\xpos,\ypos){#5}%
\fi}}
\def\putvmorphism(#1,#2)[#3`#4`#5]#6#7#8{{%
\xpos #1 \ypos #2 \arrowlength #6 \arrowtype #7
\settowidth{\xlen}{$#5$}%
\putbox(\xpos,\ypos){#3}%
{\advance \ypos by-\arrowlength
\putbox(\xpos,\ypos){#4}}%
{\advance\arrowlength by-140 \advance \ypos by-70 \ifdim\xlen>0pt
   \if m#8%
      \putsplitvector(\xpos,\ypos)\arrowlength\arrowtype
   \else
   \putvector(\xpos,\ypos)(0,-1)\arrowlength\arrowtype
   \fi
\else
   \putvector(\xpos,\ypos)(0,-1)\arrowlength\arrowtype
\fi}%
\ifdim\xlen>0pt
   \divide \arrowlength by2
   \advance\ypos by-\arrowlength
   \if l#8%
      \advance \xpos by-40
      \putrbox(\xpos,\ypos){#5}%
   \else\if r#8%
      \advance \xpos by40
      \putlbox(\xpos,\ypos){#5}%
   \else
      \putbox(\xpos,\ypos){#5}%
   \fi\fi
\fi }}
\def\putsquarep<#1>(#2)[#3;#4`#5`#6`#7]{{%
\setsqparms[#1]%
\setpos(#2)%
\settokens[#3]%
\puthmorphism(\xpos,\ypos)[\tokenc`\tokend`{#7}]{\width}{\arrowtyped}b%
\advance\ypos by \height
\puthmorphism(\xpos,\ypos)[\tokena`\tokenb`{#4}]{\width}{\arrowtypea}a%
\putvmorphism(\xpos,\ypos)[``{#5}]{\height}{\arrowtypeb}l%
\advance\xpos by \width
\putvmorphism(\xpos,\ypos)[``{#6}]{\height}{\arrowtypec}r%
}}
\def\putsquare{\@ifnextchar <{\putsquarep}{\putsquarep%
   <\arrowtypea`\arrowtypeb`\arrowtypec`\arrowtyped;\width`\height>}}
\def\square{\@ifnextchar< {\squarep}{\squarep
   <\arrowtypea`\arrowtypeb`\arrowtypec`\arrowtyped;\width`\height>}}
\def\squarep<#1>[#2`#3`#4`#5;#6`#7`#8`#9]{{
\setsqparms[#1]
\diagram
\putsquarep<\arrowtypea`\arrowtypeb`\arrowtypec`
\arrowtyped;\width`\height>
(0,0)[#2`#3`#4`{#5};#6`#7`#8`{#9}]
\enddiagram
}}                                                 
\def\putptrianglep<#1>(#2,#3)[#4`#5`#6;#7`#8`#9]{{%
\settriparms[#1]%
\xpos=#2 \ypos=#3 \advance\ypos by \height
\puthmorphism(\xpos,\ypos)[#4`#5`{#7}]{\height}{\arrowtypea}a%
\putvmorphism(\xpos,\ypos)[`#6`{#8}]{\height}{\arrowtypeb}l%
\advance\xpos by\height
\putmorphism(\xpos,\ypos)(-1,-1)[``{#9}]{\height}{\arrowtypec}r%
}}
\def\putptriangle{\@ifnextchar <{\putptrianglep}{\putptrianglep
   <\arrowtypea`\arrowtypeb`\arrowtypec;\height>}}
\def\ptriangle{\@ifnextchar <{\ptrianglep}{\ptrianglep
   <\arrowtypea`\arrowtypeb`\arrowtypec;\height>}}
\def\ptrianglep<#1>[#2`#3`#4;#5`#6`#7]{{
\settriparms[#1]
\diagram
\putptrianglep<\arrowtypea`\arrowtypeb`
\arrowtypec;\height>
(0,0)[#2`#3`#4;#5`#6`{#7}]
\enddiagram
}}                                            
\def\putqtrianglep<#1>(#2,#3)[#4`#5`#6;#7`#8`#9]{{%
\settriparms[#1]%
\xpos=#2 \ypos=#3 \advance\ypos by\height
\puthmorphism(\xpos,\ypos)[#4`#5`{#7}]{\height}{\arrowtypea}a%
\putmorphism(\xpos,\ypos)(1,-1)[``{#8}]{\height}{\arrowtypeb}l%
\advance\xpos by\height
\putvmorphism(\xpos,\ypos)[`#6`{#9}]{\height}{\arrowtypec}r%
}}
\def\putqtriangle{\@ifnextchar <{\putqtrianglep}{\putqtrianglep
   <\arrowtypea`\arrowtypeb`\arrowtypec;\height>}}
\def\qtriangle{\@ifnextchar <{\qtrianglep}{\qtrianglep
   <\arrowtypea`\arrowtypeb`\arrowtypec;\height>}}
\def\qtrianglep<#1>[#2`#3`#4;#5`#6`#7]{{
\settriparms[#1]
\width=\height                                
\diagram
\putqtrianglep<\arrowtypea`\arrowtypeb`
\arrowtypec;\height>
(0,0)[#2`#3`#4;#5`#6`{#7}]
\enddiagram
}}
\def\putdtrianglep<#1>(#2,#3)[#4`#5`#6;#7`#8`#9]{{%
\settriparms[#1]%
\xpos=#2 \ypos=#3
\puthmorphism(\xpos,\ypos)[#5`#6`{#9}]{\height}{\arrowtypec}b%
\advance\xpos by \height \advance\ypos by\height
\putmorphism(\xpos,\ypos)(-1,-1)[``{#7}]{\height}{\arrowtypea}l%
\putvmorphism(\xpos,\ypos)[#4``{#8}]{\height}{\arrowtypeb}r%
}}
\def\putdtriangle{\@ifnextchar <{\putdtrianglep}{\putdtrianglep
   <\arrowtypea`\arrowtypeb`\arrowtypec;\height>}}
\def\dtriangle{\@ifnextchar <{\dtrianglep}{\dtrianglep
   <\arrowtypea`\arrowtypeb`\arrowtypec;\height>}}
\def\dtrianglep<#1>[#2`#3`#4;#5`#6`#7]{{
\settriparms[#1]
\width=\height                                
\diagram
\putdtrianglep<\arrowtypea`\arrowtypeb`
\arrowtypec;\height>
(0,0)[#2`#3`#4;#5`#6`{#7}]
\enddiagram
}}
\def\putbtrianglep<#1>(#2,#3)[#4`#5`#6;#7`#8`#9]{{%
\settriparms[#1]%
\xpos=#2 \ypos=#3
\puthmorphism(\xpos,\ypos)[#5`#6`{#9}]{\height}{\arrowtypec}b%
\advance\ypos by\height
\putmorphism(\xpos,\ypos)(1,-1)[``{#8}]{\height}{\arrowtypeb}r%
\putvmorphism(\xpos,\ypos)[#4``{#7}]{\height}{\arrowtypea}l%
}}
\def\putbtriangle{\@ifnextchar <{\putbtrianglep}{\putbtrianglep
   <\arrowtypea`\arrowtypeb`\arrowtypec;\height>}}
\def\btriangle{\@ifnextchar <{\btrianglep}{\btrianglep
   <\arrowtypea`\arrowtypeb`\arrowtypec;\height>}}
\def\btrianglep<#1>[#2`#3`#4;#5`#6`#7]{{
\settriparms[#1]
\width=\height                               
\diagram
\putbtrianglep<\arrowtypea`\arrowtypeb`
\arrowtypec;\height>
(0,0)[#2`#3`#4;#5`#6`{#7}]
\enddiagram
}}
\def\putAtrianglep<#1>(#2,#3)[#4`#5`#6;#7`#8`#9]{{%
\settriparms[#1]%
\xpos=#2 \ypos=#3 {\multiply \height by2
\puthmorphism(\xpos,\ypos)[#5`#6`{#9}]{\height}{\arrowtypec}b}%
\advance\xpos by\height \advance\ypos by\height
\putmorphism(\xpos,\ypos)(-1,-1)[#4``{#7}]{\height}{\arrowtypea}l%
\putmorphism(\xpos,\ypos)(1,-1)[``{#8}]{\height}{\arrowtypeb}r%
}}
\def\putAtriangle{\@ifnextchar <{\putAtrianglep}{\putAtrianglep
   <\arrowtypea`\arrowtypeb`\arrowtypec;\height>}}
\def\Atriangle{\@ifnextchar <{\Atrianglep}{\Atrianglep
   <\arrowtypea`\arrowtypeb`\arrowtypec;\height>}}
\def\Atrianglep<#1>[#2`#3`#4;#5`#6`#7]{{
\settriparms[#1]
\width=\height                                     
\diagram
\putAtrianglep<\arrowtypea`\arrowtypeb`
\arrowtypec;\height>
(0,0)[#2`#3`#4;#5`#6`{#7}]
\enddiagram
}}
\def\putAtrianglepairp<#1>(#2)[#3;#4`#5`#6`#7`#8]{{%
\settripairparms[#1]%
\setpos(#2)%
\settokens[#3]%
\puthmorphism(\xpos,\ypos)[\tokenb`\tokenc`{#7}]{\height}{\arrowtyped}b%
\advance\xpos by\height
\puthmorphism(\xpos,\ypos)[\phantom{\tokenc}`\tokend`{#8}]%
{\height}{\arrowtypee}b%
\advance\ypos by\height
\putmorphism(\xpos,\ypos)(-1,-1)[\tokena``{#4}]{\height}{\arrowtypea}l%
\putvmorphism(\xpos,\ypos)[``{#5}]{\height}{\arrowtypeb}m%
\putmorphism(\xpos,\ypos)(1,-1)[``{#6}]{\height}{\arrowtypec}r%
}}
\def\putAtrianglepair{\@ifnextchar <{\putAtrianglepairp}{\putAtrianglepairp%
   <\arrowtypea`\arrowtypeb`\arrowtypec`\arrowtyped`\arrowtypee;\height>}}
\def\Atrianglepair{\@ifnextchar <{\Atrianglepairp}{\Atrianglepairp%
   <\arrowtypea`\arrowtypeb`\arrowtypec`\arrowtyped`\arrowtypee;\height>}}
\def\Atrianglepairp<#1>[#2;#3`#4`#5`#6`#7]{{
\settripairparms[#1]
\settokens[#2]
\width=\height                                
\diagram
\putAtrianglepairp                            
<\arrowtypea`\arrowtypeb`\arrowtypec`
\arrowtyped`\arrowtypee;\height>
(0,0)[{#2};#3`#4`#5`#6`{#7}]
\enddiagram
}}
\def\putVtrianglep<#1>(#2,#3)[#4`#5`#6;#7`#8`#9]{{%
\settriparms[#1]%
\xpos=#2 \ypos=#3 \advance\ypos by\height {\multiply\height by2
\puthmorphism(\xpos,\ypos)[#4`#5`{#7}]{\height}{\arrowtypea}a}%
\putmorphism(\xpos,\ypos)(1,-1)[`#6`{#8}]{\height}{\arrowtypeb}l%
\advance\xpos by\height \advance\xpos by\height
\putmorphism(\xpos,\ypos)(-1,-1)[``{#9}]{\height}{\arrowtypec}r%
}}
\def\putVtriangle{\@ifnextchar <{\putVtrianglep}{\putVtrianglep
   <\arrowtypea`\arrowtypeb`\arrowtypec;\height>}}
\def\Vtriangle{\@ifnextchar <{\Vtrianglep}{\Vtrianglep
   <\arrowtypea`\arrowtypeb`\arrowtypec;\height>}}
\def\Vtrianglep<#1>[#2`#3`#4;#5`#6`#7]{{
\settriparms[#1]
\width=\height                                 
\diagram
\putVtrianglep<\arrowtypea`\arrowtypeb`
\arrowtypec;\height>
(0,0)[#2`#3`#4;#5`#6`{#7}]
\enddiagram
}}
\def\putVtrianglepairp<#1>(#2)[#3;#4`#5`#6`#7`#8]{{
\settripairparms[#1]%
\setpos(#2)%
\settokens[#3]%
\advance\ypos by\height
\putmorphism(\xpos,\ypos)(1,-1)[`\tokend`{#6}]{\height}{\arrowtypec}l%
\puthmorphism(\xpos,\ypos)[\tokena`\tokenb`{#4}]{\height}{\arrowtypea}a%
\advance\xpos by\height
\puthmorphism(\xpos,\ypos)[\phantom{\tokenb}`\tokenc`{#5}]%
{\height}{\arrowtypeb}a%
\putvmorphism(\xpos,\ypos)[``{#7}]{\height}{\arrowtyped}m%
\advance\xpos by\height
\putmorphism(\xpos,\ypos)(-1,-1)[``{#8}]{\height}{\arrowtypee}r%
}}
\def\putVtrianglepair{\@ifnextchar <{\putVtrianglepairp}{\putVtrianglepairp%
    <\arrowtypea`\arrowtypeb`\arrowtypec`\arrowtyped`\arrowtypee;\height>}}
\def\Vtrianglepair{\@ifnextchar <{\Vtrianglepairp}{\Vtrianglepairp%
    <\arrowtypea`\arrowtypeb`\arrowtypec`\arrowtyped`\arrowtypee;\height>}}
\def\Vtrianglepairp<#1>[#2;#3`#4`#5`#6`#7]{{
\settripairparms[#1]
\settokens[#2]
\diagram
\putVtrianglepairp                             
<\arrowtypea`\arrowtypeb`\arrowtypec`
\arrowtyped`\arrowtypee;\height>
(0,0)[{#2};#3`#4`#5`#6`{#7}]
\enddiagram
}}

\def\putCtrianglep<#1>(#2,#3)[#4`#5`#6;#7`#8`#9]{{%
\settriparms[#1]%
\xpos=#2 \ypos=#3 \advance\ypos by\height
\putmorphism(\xpos,\ypos)(1,-1)[``{#9}]{\height}{\arrowtypec}l%
\advance\xpos by\height \advance\ypos by\height
\putmorphism(\xpos,\ypos)(-1,-1)[#4`#5`{#7}]{\height}{\arrowtypea}l%
{\multiply\height by 2
\putvmorphism(\xpos,\ypos)[`#6`{#8}]{\height}{\arrowtypeb}r}%
}}
\def\putCtriangle{\@ifnextchar <{\putCtrianglep}{\putCtrianglep
    <\arrowtypea`\arrowtypeb`\arrowtypec;\height>}}
\def\Ctriangle{\@ifnextchar <{\Ctrianglep}{\Ctrianglep
    <\arrowtypea`\arrowtypeb`\arrowtypec;\height>}}
\def\Ctrianglep<#1>[#2`#3`#4;#5`#6`#7]{{
\settriparms[#1]
\width=\height                               
\diagram
\putCtrianglep<\arrowtypea`\arrowtypeb`
\arrowtypec;\height>
(0,0)[#2`#3`#4;#5`#6`{#7}]
\enddiagram
}}                                           
\def\putDtrianglep<#1>(#2,#3)[#4`#5`#6;#7`#8`#9]{{%
\settriparms[#1]%
\xpos=#2 \ypos=#3 \advance\xpos by\height \advance\ypos by\height
\putmorphism(\xpos,\ypos)(-1,-1)[``{#9}]{\height}{\arrowtypec}r%
\advance\xpos by-\height \advance\ypos by\height
\putmorphism(\xpos,\ypos)(1,-1)[`#5`{#8}]{\height}{\arrowtypeb}r%
{\multiply\height by 2
\putvmorphism(\xpos,\ypos)[#4`#6`{#7}]{\height}{\arrowtypea}l}%
}}
\def\putDtriangle{\@ifnextchar <{\putDtrianglep}{\putDtrianglep
    <\arrowtypea`\arrowtypeb`\arrowtypec;\height>}}
\def\Dtriangle{\@ifnextchar <{\Dtrianglep}{\Dtrianglep
   <\arrowtypea`\arrowtypeb`\arrowtypec;\height>}}
\def\Dtrianglep<#1>[#2`#3`#4;#5`#6`#7]{{
\settriparms[#1]
\width=\height                              
\diagram
\putDtrianglep<\arrowtypea`\arrowtypeb`
\arrowtypec;\height>
(0,0)[#2`#3`#4;#5`#6`{#7}]
\enddiagram
}}                                          
\def\setrecparms[#1`#2]{\width=#1 \height=#2}%
\def\recursep<#1`#2>[#3;#4`#5`#6`#7`#8]{{%
\width=#1 \height=#2 \settokens[#3]
\settowidth{\tempdimen}{$\tokena$} \ifdim\tempdimen=0pt
  \savebox{\tempboxa}{\hbox{$\tokenb$}}%
  \savebox{\tempboxb}{\hbox{$\tokend$}}%
  \savebox{\tempboxc}{\hbox{$#6$}}%
\else
  \savebox{\tempboxa}{\hbox{$\hbox{$\tokena$}\times\hbox{$\tokenb$}$}}%
  \savebox{\tempboxb}{\hbox{$\hbox{$\tokena$}\times\hbox{$\tokend$}$}}%
  \savebox{\tempboxc}{\hbox{$\hbox{$\tokena$}\times\hbox{$#6$}$}}%
\fi \ypos=\height \divide\ypos by 2 \xpos=\ypos \advance\xpos by
\width \bfig
\putCtrianglep<-1`1`1;\ypos>(0,0)[`\tokenc`;#5`#6`{#7}]%
\puthmorphism(\ypos,0)[\tokend`\usebox{\tempboxb}`{#8}]{\width}{-1}b%
\puthmorphism(\ypos,\height)[\tokenb`\usebox{\tempboxa}`{#4}]{\width}{-1}a%
\advance\ypos by \width
\putvmorphism(\ypos,\height)[``\usebox{\tempboxc}]{\height}1r%
\efig }}
\def\recurse{\@ifnextchar <{\recursep}{\recursep<\width`\height>}}
\def\puttwohmorphisms(#1,#2)[#3`#4;#5`#6]#7#8#9{{%
\puthmorphism(#1,#2)[#3`#4`]{#7}0a \ypos=#2 \advance\ypos by 20
\puthmorphism(#1,\ypos)[\phantom{#3}`\phantom{#4}`#5]{#7}{#8}a
\advance\ypos by -40
\puthmorphism(#1,\ypos)[\phantom{#3}`\phantom{#4}`#6]{#7}{#9}b }}
\def\puttwovmorphisms(#1,#2)[#3`#4;#5`#6]#7#8#9{{%
\putvmorphism(#1,#2)[#3`#4`]{#7}0a \xpos=#1 \advance\xpos by -20
\putvmorphism(\xpos,#2)[\phantom{#3}`\phantom{#4}`#5]{#7}{#8}l
\advance\xpos by 40
\putvmorphism(\xpos,#2)[\phantom{#3}`\phantom{#4}`#6]{#7}{#9}r }}
\def\puthcoequalizer(#1)[#2`#3`#4;#5`#6`#7]#8#9{{%
\setpos(#1)%
\puttwohmorphisms(\xpos,\ypos)[#2`#3;#5`#6]{#8}11%
\advance\xpos by #8
\puthmorphism(\xpos,\ypos)[\phantom{#3}`#4`#7]{#8}1{#9} }}
\def\putvcoequalizer(#1)[#2`#3`#4;#5`#6`#7]#8#9{{%
\setpos(#1)%
\puttwovmorphisms(\xpos,\ypos)[#2`#3;#5`#6]{#8}11%
\advance\ypos by -#8
\putvmorphism(\xpos,\ypos)[\phantom{#3}`#4`#7]{#8}1{#9} }}
\def\putthreehmorphisms(#1)[#2`#3;#4`#5`#6]#7(#8)#9{{%
\setpos(#1) \settypes(#8)
\if a#9 %
     \vertsize{\tempcounta}{#5}%
     \vertsize{\tempcountb}{#6}%
     \ifnum \tempcounta<\tempcountb \tempcounta=\tempcountb \fi
\else
     \vertsize{\tempcounta}{#4}%
     \vertsize{\tempcountb}{#5}%
     \ifnum \tempcounta<\tempcountb \tempcounta=\tempcountb \fi
\fi \advance \tempcounta by 60
\puthmorphism(\xpos,\ypos)[#2`#3`#5]{#7}{\arrowtypeb}{#9}
\advance\ypos by \tempcounta
\puthmorphism(\xpos,\ypos)[\phantom{#2}`\phantom{#3}`#4]{#7}{\arrowtypea}{#9}
\advance\ypos by -\tempcounta \advance\ypos by -\tempcounta
\puthmorphism(\xpos,\ypos)[\phantom{#2}`\phantom{#3}`#6]{#7}{\arrowtypec}{#9}
}}
\def\setarrowtoks[#1`#2`#3`#4`#5`#6]{%
\def\toka{#1}
\def\tokb{#2}
\def\tokc{#3}
\def\tokd{#4}
\def\toke{#5}
\def\tokf{#6}
}
\def\hex{\@ifnextchar <{\hexp}{\hexp<1000`400>}}
\def\hexp<#1`#2>[#3`#4`#5`#6`#7`#8;#9]{%
\setarrowtoks[#9] \yext=#2 \advance \yext by #2 \xext=#1
\advance\xext by \yext \bfig
\putCtriangle<-1`0`1;#2>(0,0)[`#5`;\tokb``\tokd] \xext=#1
\yext=#2 \advance \yext by #2
\putsquare<1`0`0`1;\xext`\yext>(#2,0)[#3`#4`#7`#8;\toka```\tokf]
\advance \xext by #2
\putDtriangle<0`1`-1;#2>(\xext,0)[`#6`;`\tokc`\toke] \efig }

\newcommand{\cA}{{\mathcal A} }
\newcommand{\cB}{{\mathcal B} }
\newcommand{\cD}{{\mathcal D} }
\newcommand{\cE}{{\mathcal E} }
\newcommand{\cF}{{\mathcal F} }

\newcommand{\cH}{{\mathcal H} }
\newcommand{\cI}{{\mathcal I} }
\newcommand{\cJ}{{\mathcal J} }
\newcommand{\cK}{{\mathcal K} }
\newcommand{\cL}{{\mathcal L} }
\newcommand{\cM}{{\mathcal M} }

\newcommand{\cO}{{\mathcal O} }
\newcommand{\cP}{{\mathcal P} }

\newcommand{\cT}{{\mathcal T} }

\newcommand{\cX}{{\mathcal X} }
\newcommand{\cY}{{\mathcal Y} }

\newcommand{\cC}{{\mathcal C} }
\newcommand{\bD}{{\bf D} }

\newcommand{\wt}{\widetilde}
\newcommand{\wh}{\widehat}
\newcommand{\Hom}{{\rm Hom}}
\def\ol#1{{\overline{#1}}}
\def\psh{{plurisubharmonic}}

{\it}{\rm}
{\it}{\rm}
{\it}{\rm}
\newtheorem{corollary}{Corollary}{\bf}{\it}
\newtheorem{definition}{Definition}{\bf}{\it}
{\it}{\rm}
{\bf}{\rm}
\newtheorem{lemma}{Lemma}{\bf}{\it}
{\it}{\rm}
{\bf}{\rm}
{\it}{\rm}
\newtheorem{proposition}{Proposition}{\bf}{\it}
{\it}{\rm}
{\bf}{\rm}
\newtheorem{theorem}{Theorem}{\bf}{\it}
\newtheorem{remark}{Remark}{\it}{\rm}
\newtheorem*{condition}{Condition}{\it}{\rm}
\def\ke{K"ah\-ler-Ein\-stein }
\def\pw{Pe\-ters\-son-Weil }
\def\ks{Ko\-dai\-ra-Spen\-cer }

\begin{document}
\title[Quasi-Projectivity of Moduli Spaces]
{Quasi-Projectivity of Moduli Spaces \\ of Polarized Varieties}
\author{Georg Schumacher}
\address{Fachbereich Mathematik der Philipps-Universit"at,
Hans-Meerwein-Strasse, Lahnberge, D-35032 Marburg, Germany}
\email{schumac@mathematik.uni-marburg.de}
\author{Hajime Tsuji}
\address{Department of Mathematics,
Tokyo Institute of Technology,
2-12-1 Ohokayama, Meguro 152--8551, Japan}
\email{tsuji@math.titech.ac.jp}
\begin{abstract}
By means of analytic methods the quasi-projectivity of the moduli
space of algebraically polarized varieties with a not necessarily
reduced complex structure is proven including the case of non-uniruled
polarized varieties.

\end{abstract}
\maketitle
\tableofcontents

\section{Introduction} In the theory of algebraic varieties, it is
fundamental to study the moduli spaces of algebraic varieties. As
for the existence of moduli spaces, it had been known that there
exists an algebraic space as a coarse moduli space of
non-uniruled polarized projective manifolds with a given Hilbert
polynomial. Here an algebraic space denotes a space which is
locally a finite quotient of an algebraic variety. Actually the
notion of algebraic spaces was introduced to describe the moduli
spaces (\cite{ar}). According to the theory of algebraic spaces
by M.~Artin (\cite{ar, ar2, knt}), the category of proper
algebraic spaces of finite type defined over $\mathbb C$ is
equivalent to the category of Moishezon spaces. Hence the moduli
spaces of non-uniruled polarized manifolds have abundant
meromorphic functions and were considered to be not far from being
quasiprojective. Various attempts were made to prove the
quasiprojectivity of the moduli spaces of non-uniruled, polarized
algebraic varieties (cf.\ \cite{k-m, kn, ko, v}). E.~Viehweg
(\cite{v}) developed a theory to construct positive line bundles
on moduli spaces. He used results on the weak semipositivity of
the direct images of relative multicanonical bundles. In
particular he could prove the quasiprojectivity of the moduli
spaces of canonically polarized manifolds (\cite{v}).
J.~Koll\'{a}r studied the Nakai-Moishezon criterion for ampleness
on certain complete moduli spaces in \cite{ko}, with applications
to the projectivity of the moduli space of stable curves and
certain moduli spaces of stable surfaces under boundedness
conditions. However, his approach appears quite different from
our present methods, which do not require the completeness of
moduli spaces. His result was used to show the projectivity of
the compactified moduli spaces of surfaces with ample canonical
bundles by V.~Alexeev (\cite{al}).

The main result in this paper is the quasi-projectivity of the
moduli space of non-uniruled polarized manifolds. Non-uniruledness
is in fact not used here. All we need is the existence of a
moduli space.

In fact, given a polarized projective manifold, a universal
family of embedded projective manifolds over a Zariski open
subspace $\cH$ of a Hilbert scheme is determined after fixing the
Hilbert polynomial.

The identification of points of $\cH$, whose fibers are
isomorphic as polarized varieties, defines an analytic equivalence
relation $\sim$ such that the set theoretic moduli space is
$\cM=\cH/\!\sim$. The quotient is already a complex space, if the
equivalence relation is proper. Moreover,in this situation, it
follows that $\cM$ is an algebraic space. If the above equivalence
relation is induced by the action of a projective linear group
$G$, properness of $\sim$ means properness of the action of $G$.
In this moduli theoretic case $\cH/\!\sim$ is already a geometric
quotient.

\begin{theorem}\label{main}
Let $ \cK$ be a class of polarized, projective manifolds such
that the moduli space $\cM$ exists as a proper quotient of a
Zariski open subspace of a Hilbert scheme. Then $\cM$ is
quasi-projective.
\end{theorem}

The proof of the theorem consists of two steps. The first step is
to construct a line bundle on the compactified moduli space with a
singular hermitian metric of strictly positive curvature on the
interior.

The method is based upon the curvature formula for Quillen metrics
on determinant line bundles (\cite{bgs}), the theory of Griffiths
about period mappings (\cite{griff}), and moduli of framed manifolds.

The second step is to construct sufficiently many holomorphic
sections of a power of the above line bundles in terms of
$L^{2}$-estimates of the $\ol{\partial}$-operator. The key
ingredient here is the theory of closed positive
$(1,1)$-currents, which controls the multiplier ideal sheaf of a
singular hermitian metric. This step can be viewed as an
extension of the Kodaira embedding theorem to the quasi-projective
case.

\section{Singular hermitian metrics}\label{singherm}

\begin{definition} Let $X$ be a complex manifold and $L$ a
holomorphic line bundle on $X$. Let $h_0$ be a hermitian metric
on $L$ of class $C^\infty$ and $\varphi \in L^1_{loc}(X)$. Then
$h = h_0 \cdot {\bf e}^{-\varphi}$ is called a singular hermitian
metric on $L$.
\end{definition}

Following the notation of \cite{demai4} we set
$$
d^c= \frac{\sqrt{-1}}{2\pi}(\partial - \overline{\partial})
$$
and call the real $(1,1)$-current
\begin{equation}\label{poscurr}
\Theta_h = d d^c (-\log h) = -\frac{\sqrt{-1}}{\pi}\partial
\ol\partial \log h
\end{equation}
the ''curvature current'' of $h$. It differs from the Chern
current by a factor of $2$.

A real current $\Theta$ of type $(1,1)$ on a complex manifold of
dimension $n$ is called {\em positive}, if for all smooth
$(1,0)$-forms $\alpha_2, \ldots, \alpha_n$
$$
\Theta \wedge \sqrt{-1}\alpha_2 \wedge \ol \alpha_2 \wedge \ldots
\wedge \sqrt{-1}\alpha_n \wedge \ol \alpha_n
$$
is a positive measure. We write $\Theta \geq 0$.

A singular hermitian metric $h$ with positive curvature current is
called {\em positive}. This condition is equivalent to saying
that the locally defined function $-\log h$ is \psh.

Let $W\subset \mathbb C^n$ be a domain, and $\Theta$ a positive
current of degree $(q,q)$ on $W$. For a point $p\in W$ one defines
$$
\nu(\Theta,p,r)=\frac{1}{r^{2(n-q)}}\int_{\|z-p\|<r} \Theta(z)
\wedge (dd^c\|z\|^2)^{n-q}.
$$
The {\it Lelong number} of $\Theta$ at $p$ is defined as
$$
\nu(\Theta,p)= \lim_{{r\to 0}\atop{r\; >\; 0}}\nu(\Theta,p,r).
$$
If $\Theta$ is the curvature of $h={\bf e}^{-u}$, $u$
plurisubharmonic, one has
$$
\nu(\Theta,p)= \sup \{\gamma \geq 0; u\leq\gamma\log(\|z-p\|^2) +
O(1)\}.
$$

The definition of a singular hermitian metric carries over to the
situation of reduced complex spaces.

\begin{definition}\label{singdef} Let $Z$ be a reduced complex space
and $L$ a holomorphic line bundle. A {\em singular} hermitian
metric $h$ on $L$ is a singular hermitian metric $h$ on
$L|Z_{reg}$ with the following property: There exists a
desingularization $\pi: \wt{Z} \to Z$ such that $h$ can be
extended from $Z_{reg}$ to a singular hermitian metric $\wt{h}$
on $\pi^*L$ over $\wt{Z}$.
\end{definition}

The definition is independent of the choice of a desingularization
under a further assumption. Suppose that $\Theta_{\wt{h}} \geq -c
\cdot \omega $ in the sense of currents, where $c>0$, and $\omega$
is a positive definite, real $(1,1)$-form on $\wt{Z}$ of class
$C^{\infty}$. Let $ \pi_1 :{Z_1} \to Z$ be a further
desingularization. Then $\wt{Z}\times_Z Z_1 \to Z$ is dominated
by a desingularization $Z'$ with projections $p :Z' \to \wt{Z} $
and $p_1 :Z' \to Z_1$. Now $p^*\log \wt{h}$ is of class
$L^1_{loc}$ on $Z'$ with a similar lower estimate for the
curvature. The push-forward $p_{1*}p^*\wt{h}$ is a singular
hermitian metric on $Z_1$. In particular, the extension of $h$ to
a desingularization of $Z$ is unique. \qed

In \cite{g-r} for plurisubharmonic functions on a normal complex
space the Riemann extension theorems were proved, which will be
essential for our application. The relationship with the theory
of distributions was treated in \cite{demai}.

For a reduced complex space a {\it \psh} function $u$ is by
definition an upper semi-continuous function $u:X \to [-\infty,
\infty)$ whose restriction to any local, smoothly parameterized
analytic curve is either identically $-\infty$ or subharmonic.

A locally bounded function $u:X \to [-\infty, \infty)$ from
$L^1_{loc}(X)$ is called {\it weakly \psh}, if its restriction to
the regular part of $X$ is \psh.

Differential forms with compact support on a reduced complex
space are by definition locally extendable to an ambient
subspace, which is an open subset $U$ of some $\mathbb C^n$. Hence
the dual spaces of differential $C^\infty$-forms on such $U$
define currents on analytic subsets of $U$. The positivity of a
real $(1,1)$-current is defined in a similar way as above
involving expressions of the form (\ref{poscurr}).

For locally bounded functions of class $L^1_{loc}$ the weak
plurisubharmonicity is equivalent to the positivity of the
current $d d^c u$. It was shown that these functions are exactly
those, whose pull-back to the normalization of $X$ are \psh. We
note

\begin{definition}
Let $L$ be a holomorphic line bundle on a reduced complex space
$X$. Then a singular hermitian metric $h$ is called positive, if
the functions, which define locally $-\log h$ are weakly \psh.
\end{definition}

This definition is compatible with Definition~\ref{singdef}: Let
$L$ be a holomorphic line bundle on a complex space $Z$ equipped
with a positive, singular hermitian metric $h_r$ on $L|Z_{reg}$.
If $\pi:\wt Z \to Z$ is a desingularization, and $\wt h$ a
positive, singular hermitian metric on $\pi^* L$, extending
$h|Z_{reg}$, we see that $-\log h_r$ is locally bounded at the
singularities of $Z$ so that $\wt h$ induces a singular, positive
metric on $L$ over $Z$.

\section{Deformation theory of framed manifolds--
$V$-structures}\label{frV} Let $X$ be a compact complex
manifold and $D \subset X$ a smooth (irreducible) divisor. Then
$(X, D)$ is called a logarithmic pair or a {\it framed manifold}.

For {\it any} $m \in \mathbb N$ an associated $V$-structure $\wt
X_m$ on $X$ is defined, in terms of local charts $\pi: W \to U$,
$U \subset X$, $W \subset \mathbb C^n$ such that $\pi$ is just an
isomorphism, if $U \cap D = \emptyset$ or a Galois covering of
order $m$ with branch locus $U \cap D$ (and group $G=\mathbb
Z_m$).

By definition, the differential forms and vector fields on $X$
with respect to the $V$-structure, which are $V$-differentiable or
$V$-holomorphic, are defined on $X\backslash D$ with the property
that the local lifts under $\pi| W \cap \pi^{-1}(U\backslash D) :
W \cap \pi^{-1}(U\backslash D) \to U \backslash D$ can be
extended in a holomorphic or differentiable way to $W$.

With $m$ being fixed, we denote by ${\mathcal T}_X^V$ and
${\mathcal A}^{V, q}_X({\mathcal T}^V_X)$ resp.\ the sheaves of
$V$-holomorphic vector fields and $V$-differentiable $q$-forms
with values in ${\mathcal T}_X^V$ resp.
\begin{lemma}\label{dol}
\begin{enumerate}
\item[(i)] For any $m\in \mathbb N$ the Dolbeault complex
$$
0 \to {\mathcal T}^V_X \to {\mathcal A}^{V, \bullet}_X({\mathcal
T}^V_X)
$$
is well-defined and exact.
\item[(ii)] The sheaf ${\mathcal T}_X^V$ is canonically isomorphic to
$\Omega^1_X(\log D)^{\wedge}$.
\end{enumerate}
\end{lemma}

By definition, a family $(\cX_s,\cD_s)_{s\in S}$ of framed
manifolds, parameterized by a complex space $S$ is given by a
smooth, proper, holomorphic map $f:\cX \to S$ together with a
divisor $\cD \subset \cX$, such that $f|\cD$ is proper and smooth,
such that $\cX_s= f^{-1}(s)$ and $\cD_s= \cD \cap \cX_s$. A local
deformation of a framed manifold $(X,D)$ over a complex space $S$
with base point $s_0\in S$ is a deformation of the embedding $i:
D \hookrightarrow X$, i.e. induced by a family $\cD
\hookrightarrow \cX \to S$ together with an isomorphism
$(X,D)\stackrel{\sim}{\to} (\cX_{s_0}, \cD_{s_0})$, where two
such objects are identified, if these are isomorphic over a
neighborhood of the base point. The existence of versal
deformations (i.e.\ complete and semi-universal deformations) of
these objects is known. We denote by $T^\bullet(X)\simeq
H^\bullet(X,\cT_X)$ and $T^\bullet(X,D)$ resp.\ the tangent
cohomology of $X$ and $(X,D)$ resp.
\begin{corollary}
The space of infinitesimal deformations of $(X,D)$ equals
$
T^1(X,D)=H^1(\Gamma(X,{\mathcal A}^{V, \bullet}_X({\mathcal
T}^V_X))).
$
It can also be computed in terms of Cech cohomology as
$\overset{\vee}{H^1}(\mathfrak{U}, \cT_X^V)$ of $V$-holomorphic
vector fields, where $\mathfrak U$ is a $G$-invariant
$\cT^V_X$-acyclic covering.
\end{corollary}
We have the following exact sequence.:
\begin{gather}\label{longtangcoh}
0 \to T^0(X,D) \to T^0(X) \to H^0(D,\cO_D(D)) \to \\
\hspace{2cm} T^1(X,D) \to T^1(X) \to H^1(D,\cO_D(D))\notag.
\end{gather}
We denote by $T^1_1(X)\subset T^1(X)$ the image of $T^1(X,D)$. The
composition of $H^1(X,\cT_X)) \to H^1(D,\cO_D(D))$ with the
natural map $H^1(D,\cO_D(D)\to H^2(X,\cO_X)$ equals the map
induced by the cup-product with the Chern class of $D$. The latter
is induced by the Atiyah sequence for the pair $(X,\cO_X(D))$, and
its kernel $T^1_0(X)$ consists of those infinitesimal deformations
for which the isomorphism class of the line bundle $[D]$ extends.
Assume that  $D$ is an ample divisor on $X$, and
$\lambda_X=c_1(D)$ its (real) Chern class. Then the pair
$(X,\lambda_X)$ is a {\it polarized variety}, and $T_0^1(X)$ is
the space of infinitesimal deformations of $(X,\lambda_X)$.
Studying moduli spaces of polarized varieties, we are free to
replace the ample divisor $D$ by a uniformly chosen multiple, in
which case $T_0^1(X)$ and $T_1^1(X)$ can be identified.

\begin{remark} The group of infinitesimal automorphisms $T^0(X, D)$
vanishes, if $K_X + [D]$ is positive. Like in the case of
canonically polarized manifolds, in a family of such {\it framed}
manifolds the relative automorphism functor (or more generally
isomorphism functor) is represented by a space such that the
natural map to the base is finite and proper. Moreover, general
deformation theory implies that any versal deformation is
universal.
\end{remark}

\section{Cyclic coverings}\label{cycov} Let $X$ be a compact
complex manifold, and $D,D'$ effective divisors such that $D \sim
m \cdot D'$ for some $ m \in \mathbb N$. Denote by $E$ and $E'$
resp.\ bundle spaces for the corresponding line bundles. Let
\begin{equation}
\settriparms[1`1`1;300] \Vtriangle[E'`E`X;\ell`\pi'`\pi]
\end{equation}
be the morphism over $X$, which sends a bundle coordinate
$\alpha$ to $\alpha^m$.

Let $\sigma$ be a canonical section of $\pi$.  Then we define
$X_m=V(\ell - \sigma \circ \pi') \subset E'$. If $D$ is a smooth
divisor, the subspace $X_m \subset E$ is a manifold, and
$\pi'|X_m: X_m \to X$ is a Galois covering with branch locus
$D\subset X$.

We assume now that $D$ is very ample, providing an embedding
$\Phi : X \to \mathbb P_N$. We denote by $P$ the dual projective
space, and by $\Sigma \subset \mathbb P_N \times P$ the
tautological hypersurface with divisor $\cD = \Sigma \cap (X
\times P) \subset X \times P$ and bundle space $\cE \to X \times
P$. Let $\cD_t = \Sigma_t \cap X$ for $t\in P$.

We have flat families over $X\times P$ and $P$ resp.
\begin{equation}\label{gal}
\begin{picture}(1000,800)
\settriparms[1`1`1;350] \putqtriangle(100,400)[{\underline
X_m}`\underline E'`X\times P;`\mu`] \settriparms[1`0`1;350]
\putptriangle(450,400)[\phantom{\underline E'}`\underline
E`\phantom{X \times P};``] \put(450,300){\vector(0,-1){200}}
\put(400,0){$P$} \put(70,660){\vector(1,-2){300}}
\put(110,350){$\pi$}
\end{picture}
\end{equation}
Here the bundle $\underline E$ comes from the globally defined
divisor $\cD$. The bundle $\underline E'$ is first defined locally
with respect to $P$, which gives rise to a cocycle with
coefficients from $\mathbb C^*$, i.e.\ a line bundle on $P$.

\begin{proposition}\label{galois} The total space $\underline X_m$ is
smooth. In particular, the dualizing sheaf
$\omega_{\underline{X}_m/P}$ equals the relative canonical sheaf
$K_{\underline{X}_m/P}:=K_{\underline{X}_m}\otimes \pi^*K^{-1}_P$.
\end{proposition}

\begin{proof}
As $\underline X_m \subset \underline{E}'$ is of codimension one,
it is sufficient, to find  a local function for any $x_0 \in
\underline X_m$, which vanishes at $x_0$, and whose gradient at
this point is non-zero. Let again $\sigma$ be a canonical section
of the line bundle $\underline E$ over $X \times P$. We denote by
$t_0$ the image of $x_0$ in $P$, and take local coordinates $t$ of
$P$ around $t_0$. Let $\alpha$ be a local bundle coordinate of
$\underline E'$ around $t_0$, and $z$ a local coordinate on $X$
so that $x_0$ is given by $(z_0,\alpha_0, t_0)$. Now $t_0\in P$
corresponds to a section $\sigma_{t_0}(z)$ of $\underline{E}|X
\times \{t_0\}$. The space $\underline X_m$ is defined by
$g(z,\alpha,t):=\sigma_t(x) - \alpha^m=0$ around $x_0$. If
$\alpha_0\neq 0$, we have $(\partial g/\partial \alpha)(x_0)\neq
0$. If $\alpha_0=0$ holds, $\sigma_{t_0}(z_0)=0$. Since $D$ is
very ample on $X$,  we find a section of $\underline E|X \times
\{t_0\}$, which does not vanish at $x_0$. This section gives some
$t_1\in P$, i.e.\ some $\sigma_{t_1}$. Let $\sigma_{t(\tau)} =
\sigma_{t_0} + \tau \sigma_{t_1}$ be the line through $t_0$ and
$t_1$. Then $(\partial g/\partial \tau)|_{\tau=0}\neq 0$.
\end{proof}
\begin{remark}
The analogous statement is true for smooth families $f:\cX \to S$
of manifolds over complex spaces, assuming the existence of
families of divisors $\cD$ and $\cD'$, with $\cD \sim m\cdot\cD'$,
and a relative embedding $\cX \hookrightarrow {\mathbb P}_N
\times S$ induced by $\cD$. We have
\begin{equation}\label{relgal}
\begin{picture}(1000,800)
\settriparms[1`1`1;350]
\putqtriangle(100,400)[{\cX_m}`\cE'`\cX\times P;`\mu`]
\settriparms[1`0`1;350]
\putptriangle(450,400)[\phantom{\cE'}`\cE`\phantom{\cX \times
P};``] \put(450,300){\vector(0,-1){200}} \put(400,0){$S\times P$}
\put(70,660){\vector(1,-2){300}} \put(130,300){$f_m$}
\put(500,200){$f\times \rm{id}$}
\end{picture}
\end{equation}
and the induced map $\cX_m \to S$ is smooth. In particular, the
canonical and dualizing sheaves $K_{\cX_m/S\times P}=
K_{\cX_m}\otimes f_m^* K_{S\times P}^{-1}$ and
$\omega_{\cX_m/S\times P}$ resp.\ are isomorphic, if $S$ is
smooth.
\end{remark}

Let $(X, D)$ be a framed manifold, and $D\sim mD'$ for some
effective $D'$ as above. Again, let $G=\mathbb Z_m$ denote the
Galois group, $X$ is isomorphic to the quotient $X_m/G$, and the
group $G$ acts on $H^1(X_m,\cT_{X_m})$ with invariant subgroup
$
H^1(X_m,\cT_{X_m}) \supset H^1(X_m,\cT_{X_m})^G.
$
The average over the group defines a retraction. Next, we identify
$H^1(X_m,\cT_{X_m})^G$ with the $V$-tangent cohomology group
$\overset{\vee}{H^1}(\mathfrak{U}, \cT_X^V)$ in the sense of
Section~\ref{frV}: The morphisms
$
C^\bullet(\mathfrak{U},\cT_{X_m})^G \hookrightarrow
C^\bullet(\mathfrak{U},\cT_{X_m})   \overset{r}{\to }
C^\bullet(\mathfrak{U},\cT_{X_m})^G
$
descend to the cohomology and $
C^\bullet(\mathfrak{U},\cT_{X_m})^G \simeq
C^\bullet(\mathfrak{U},\cT_X^V). $ This argument avoids any
smoothing of invariant differential forms.

\begin{remark}\label{infdef}
The infinitesimal deformations of a framed manifold $(X,D)$ can be
identified with
$
T^1(X,D)= H^1(\Gamma(X,{\mathcal A}^{V, \bullet}_X({\mathcal
T}^V_X))) =H^1(\mathfrak{U},\cT_X^V)=H^1(X_m,\cT_{X_m})^G.
$
\end{remark}

\section{Canonically polarized framed manifolds}\label{framf}
\begin{definition} A framed manifold $(X,D)$ is called
\begin{itemize}
\item[(i)]
 {\em canonically polarized}, if
$$
K_X +[D]>0
$$
\item[(ii)] {\em $m$-framed} under the condition
$$
(*)_m \qquad K_X + \frac{m-1}{m}[D] >0
$$
for some $m \geq 2$.
\end{itemize}
\end{definition}
In the sequel we always assume condition $(*)_m$ for some fixed
$m$. We note that for the Galois covering $\mu: X_m \to X$ with
smooth $X_m$ the relation
$$
\mu^*(K_X + \frac{m-1}{m}[D])=K_{X_m}
$$
holds. In our applications the divisor $D$ will always be ample so
that condition (ii) is slightly stronger than (i). We will still
use the term ''canonically polarized framed manifold'' in this
case. This will also be justified later.

\begin{proposition}\label{X_m}
Let $D' \subset X$ be a very ample divisor as above, and $m>2$.
Let $D \subset X$ be a smooth divisor $D\sim m\cdot D'$ such that
$$
K_X + \frac{m-2}{m} D
$$
is very ample. Then the canonical bundle $K_{X_m}$ is very ample.
\end{proposition}

\begin{proof} The sheaf $\cO_X(K_X + \frac{m-1}{m}D) \subset
\mu_*(\cO_{X_m}(K_{X_{m}}))$ is a direct summand. Let $\mathbb Z_m
\simeq G \hookrightarrow {\rm Aut}(X_m)$ be the group of deck
transformations with a generator $\gamma$, and denote by $\zeta$ a
primitive m-the root of unity. Let $\oplus_{j=1}^m E_j$ be an
eigenspace decomposition of the space of global sections of
$K_{X_m}$ with respect to the eigenvalues $\zeta^j$ of $\gamma$.
It follows that the spaces $E_j$ can be identified with the space
of global sections of $K_X + (m-j)\cdot D'$, again with $j=1,
\ldots, m$. The pull-backs of sections of such a space are
sections of $K_{X_m} - (j-1)A$, where $A \subset X_m$, $A \simeq
D'$, is the branching divisor of $\mu$, so that the
identification $\Gamma(X, \cO_X(K_X+ (m-j)\cdot D')) \simeq E_j$
is a multiplication with a canonical section of $[(j-1)A]$.

The space $E_1$ clearly separates points, whose image under $\mu$
are different.

Let $p, q \in X_m$ with $\mu(p)=\mu(q)=x$. Then there exist
sections of $[K_X + (m-2) D']$ and $[K_X + (m-1) D']$ which do
not vanish at $x$. A suitable linear combination of the induced
elements of $E_1$ and $E_2$ separates $p$ and $q$. The argument
is also applicable to tangent vectors.
\end{proof}

Now we consider the situation given in diagram~(\ref{relgal}),
where $S$ need not be smooth. Let $\cA \subset S \times P$ be the
locus of singular divisors $D$. Over its complement the direct
image of the relative canonical sheaf is certainly locally free.

We write $\cX_m':=\cX_m \backslash f_m^{-1}(\cA)$, $T:=P\times
S$, $T':=T \backslash \cA$, and $f_m'$ for the restriction of the
map $f_m$. In a similar way we restrict $\wt f:=f \times{\rm id}$
to $T'$ and get $\wt{f}':(\cX\times P)'\to T'$.

\begin{proposition}\label{exK}
The locally free sheaf $f_{m*}'K_{\cX_m'/T'}$ possesses a
natural, locally free extension.
\end{proposition}
\begin{proof}
We use the decomposition $f_{m*}K_{\cX_m'/T'}=\oplus_{j=0}^{m-1}
\wt f'_*(K_{(\cX\times P)' /T'} + j\cdot [\cD'|(\cX\times P)'])$
from the proof of Proposition~\ref{X_m}. Now for the family
$(\cX\times P)' \to T'$, with relatively (very) ample divisor
$\cD'$, the Kodaira-Nakano vanishing theorem and the
Grothendieck-Grauert comparison theorem show that for $j>0$ the
sheaves $\wt f_*(K_{(\cX\times P)/T}+ j \cdot [\cD'])$ are locally
free on $T$ (Here the divisor $\cD'$ corresponds to the line
bundle $\cE'$). Let $j=0$. Since $f'_{m*}(K_{\cX_m'/T'})$ is
locally free on $T'$, also $\wt f_*(K_{\cX\times P/T})$ is locally
free, when restricted to $T'$. On the other hand, it does not
involve the divisor $\cD'$, and $\wt f_*(K_{\cX \times P/T})$ is
the pull-back of the direct image of $K_{\cX/S}$, so it is
constant along all fibers of $T \to S$, and locally free in the
interior, hence also $\wt f_*(K_{K_{\cX\times P/T}})$ is locally
free.
\end{proof}

Next, we want to recover the above extension of the relative
canonical sheaf. We have the diagram (\ref{relgal}). The fibers of
$f_m$ are branched along the $\cD_s$ with singularities over the
singularities of the branching divisors. By definition the map
$f_m$ is flat with Cohen-Macaulay fibers. According to results of
Kleiman \cite{kleim} for such morphisms relative dualizing
sheaves commute with base change. Again, we denote by the letter
$\omega$ dualizing sheaves.

Now
\begin{equation*}
\begin{split}
f_{m*}(\omega_{\cX_m/T}) &\simeq {\it Hom}_{\cO_{T }} (R^n
f_{m*}(\cO_{\cX_m}),\cO_{T})\\ &= {\it Hom}_{\cO_{T}}((R^n
f_{*})(\mu_*
\cO_{\cX_m}),\cO_{T})\\
 &\simeq
{\it Hom}_{\cO_{T}}((R^n f_{*})(\oplus_{j=0}^{m-1}(\cO_{\cX\times
P}(-j\cdot \cD'))
,\cO_{T})\\
& \simeq
 f_*{\it Hom}_{\cO_{\cX \times P}}( \mu_*\cO_{\cX_m} ,
 \omega_{\cX\times P/T} )
\end{split}
\end{equation*}
Altogether, we have
\begin{lemma}\label{relcan}
$$
 f_{m*}(\omega_{\cX_m/T})\simeq
\oplus_{j=0}^{m-1} f_*({\omega_{\cX\times P/T}}(j\cdot \cD')).
$$
\end{lemma}
In particular, the extension  from Proposition~\ref{relcan}
equals $ f_{m*}(\omega_{\cX_m/T})$ (which is compatible with
further pull-backs).

Later we will consider this sheaf from a Hodge theoretic
viewpoint.

\section{Singular Hermitian metrics for families of canonically
polarized framed manifolds} We first recall some facts concerning
the period map in the sense of Griffiths \cite{griff} for families
$f:\cY \to S$ of manifolds with very ample canonical bundle. We
will apply the results to families of the form $f_m: \cX_m \to S$
with relative dimension $n$ from Section~\ref{cycov}. The direct
image under $f$ of the relative canonical sheaf $K_{{\cY}/S}$ is
also called Hodge bundle ${\mathcal E}_0$. It is equipped with
the flat metric from ${\rm R}^nf_*\mathbb C$. Explicitly, for any
two holomorphic $n$-forms $\phi$ and $\psi$ on a manifold
$\cY_{s_0}$ we have
$$
(\phi,\psi):= (\sqrt{-1})^{n^2}\int_{\cY_{s_0}} \phi\wedge
\overline{\psi}.
$$
Let $\partial/\partial s$ be a tangent vector at a point $s_0$
Then the contraction with the Kodaira-Spencer class
$[A^\alpha_{s_0\bar\beta} \frac{\partial}{\partial
z^\alpha}dz^{\bar\beta}] \in H^1(\cY_{s_0},\cT_{\cY_{s_0}})$
induces a linear map
\begin{equation}\label{sigma0}
\sigma_0(\frac{\partial}{\partial s}\left.\right|_{s_0}):
H^0(\cY_{s_0}, \Omega^{n}_{\cY_{s_0}}) \to H^1(\cY_{s_0},
\Omega^{n-1}_{\cY_{s_0}}).
\end{equation}
The natural metric on the latter space is again induced by the
integration of exterior products of differential forms, after we
provide the fibers with a family of auxiliary K"ahler structures
(e.g.\ of \ke type). Following Griffiths
\cite[Theorem~(5.2)]{griff} the curvature $\Theta_0$ of this
hermitian metric is given by the formula
\begin{equation}\label{curvsig}
(\Theta_{0}\phi,\psi)=  (\sqrt{-1})^{n^2} \int_{\cY_{s_0}}
H\left(\sigma_0(\frac{\partial}{\partial s})(\phi) \right)\wedge
\overline{H\left(\sigma_0(\frac{\partial}{\partial
s})(\psi)\right)},
\end{equation}
which is defined in terms of cohomology classes. (Here $H$ denotes
the harmonic projection). So $\Theta_{0}$ is semi-positive as well
as its trace ${\rm tr}(\Theta_0)$. If ${\rm tr}(\Theta_{0})
(\frac{\partial}{\partial s},\frac{\partial}{\partial
s})|_{s_0}=0$, also  $\Theta_{0} (\frac{\partial}{\partial
s},\frac{\partial}{\partial s})|_{s_0}$ vanishes. The auxiliary
K"ahler metric is only needed to show the positivity of the
curvature, the metric on the relative canonical bundle is
independent of the choice. The sheaf $R^1f_*\Omega^{n-1}_{\cY/S}$
is usually denoted by $\cE_1$.

Denote by $\bD$ the period domain of Hodge structures, and by
$\Phi : S \to \bD$ the induced (multivalued) period map. Then
$\Hom( \cE^0\otimes_{\cO_S}\mathbb C(s),
\cE^1\otimes_{\cO_S}\mathbb C(s))$ is a subspace of the tangent
space of $\bD$ at the point $\Phi(s)$, and it carries the natural
$L^2$-inner product (cf.\ (\ref{curvsig})). We denote this metric
by $ds_0^*$. If $S \simeq \Delta^{*k}\times \Delta^\ell$, one
knows $\Phi^*ds_0^* \leq {\rm const.\;} ds^2_{Poinc}$. On the
other hand, for $f:\cY \to S$, by (\ref{curvsig}), the trace of
the curvature of the flat metric restricted to a bundle $\cE_0$
gives exactly $ds^2_0$. This argument shows:

\begin{lemma}\label{delta} Let $\cY \to S$, $S=\Delta^{*k}\times
\Delta^\ell$ be a holomorphic family of canonically polarized
manifolds. Let $h_S$ be the natural $C^\infty$ hermitian metric on
$\det f_*\cK_{\cY/S}$. Then the curvature $\Theta_S$ is
semi-positive (in the sense of $C^\infty$-forms), and dominated by
a constant multiple of the K"ahler form $\omega_S$ induced by
$ds^2_{Poinc}$.
\end{lemma}

For effectively parameterized families $f_m: \cX_m \to T$ and
large $m$ the map $ \sigma_0: H^1(\cX_{m,s_0},\cT_{\cX_{m,s_0}})
\to \Hom(H^0(\cX_{s_0}, \Omega^{n}_{\cX_{s_0}}),H^1(\cX_{s_0},
\Omega^{n-1}_{\cX_{s_0}})) $ is in fact injective. This was shown
in a general setting by Ivinskis, and attributed to Griffiths in
\cite{ivin} for the special case of cycling coverings. (The
cohomological assumption for \cite[Theorem~2.4]{ivin} is easily
satisfied for sufficiently high powers $m$ of the line bundles
$[\cD_s]$, which can be chosen uniformly on the moduli space from
the beginning by Noether induction.)

Now the base $S$ is equipped with the line bundle $\lambda_{fr}=
\det f_{m*}\cK_{{\cX_m}/T}$ (which equals the determinant line
bundle in the sense of the derived category, because of the
Kodaira vanishing theorem). Then the curvature of the induced
hermitian metric $h$ on $\lambda$ is $\Theta_h={\rm
tr}(\Theta_0)$. Altogether:

\begin{proposition}\label{curvper} The curvature $\Theta_h$ of
$(\lambda_{fr}, h)$ is semi-positive. It is strictly positive in
all directions, where the family is effectively parameterized
(i.e.\ where the Kodaira-Spencer map is non-zero).
\end{proposition}

Now we return to the notation of Section~\ref{framf}. The main
theorem is stated for non-singular base spaces.

\begin{theorem}\label{lelong}
The determinant (invertible) sheaf $\det f_{m*}K_{\cX_m/T}$
carries a natural positive hermitian metric, whose Lelong numbers
vanish everywhere. Moreover, for all $p \in \mathbb N$, the
exterior powers $\Theta_h^p$ of its curvature form $\Theta_h$ are
well-defined $(p,p)$-currents, whose Lelong numbers vanish
everywhere as well.
\end{theorem}

We shall apply the theorem in two different situations: Over the
interior of the moduli space we deal with families of manifolds
of the type $X_m$, where in the limit we have singular Galois
coverings $X_m \to X$ (cf.\ Section~\ref{framf}). Here the key
point is that the total space $\cX_m$ is already smooth according
to Proposition~\ref{galois} so that we can identify the relative
dualizing sheaf with the relative canonical sheaf. The other
situation occurs at the boundary of the moduli space, where we are
free to modify the boundary.

The theorem follows from the known results in the theory of mixed
Hodge structures. We show here an upper estimate for a singular
Hermitian metric. Together with the positivity of this metric the
vanishing of the Lelong numbers follows.

Concerning singular base spaces of holomorphic families, we
observe that the $L^2$-inner products (for tangent vectors of the
base) are well-defined for singular bases spaces. For our
applications we will need that the construction is functorial,
i.e.\ compatible with base changes like restrictions to closed
subspaces and desingularizations in the view of
Definition~\ref{singdef}.

For a family $f_m:\cX_m \to T$ ($T$ is smooth), we denote by $\cA
\subset T$ the set of points with singular fibers. Let $\nu: \wt
T \to T$ be given by a sequence of blow-ups with regular centers
so that the preimage $\cB$ of $\cA$ is a normal crossings
divisor. Let $\wt \cX_m \to \cX_m \times_T \wt T$ be a
desingularization of the of the component of $\cX_m \times_T \wt
T$ that dominates $\wt T$, with the property that the preimage of
$\cB$ is a normal crossing divisor. Let
\begin{equation}
\begin{picture}(500,500)
\setsqparms[1`1`1`1;400`400]
\putsquare(100,100)[\wt\cX_m`\cX_m`\wt T`T;\nu'`\wt f_m`f_m`\nu]
\end{picture} \qquad
\end{equation}
be the induced commutative diagram. We denote by a prime accent
the restriction of fiber spaces to the resp.\ complements of
normal crossing divisors.

An argument of Deligne shows that the local monodromy of
$R^nf_{m*}\mathbb C$ on $T'$ is unipotent around generic points
of $\cA$, i.e.\ in codimension one. And since it is locally
abelian on $\wt T'$, this holds everywhere. For our purpose the
unipotent reduction is sufficient. We need a local statement with
respect to the base $T$. The argument is known: Around each
component of the normal crossings divisor $\cB$ the eigenvalues
of the local monodromy transformation on $R^n \wt f'_{m*}\mathbb
C$ are certain roots of unity \cite{borel}. After taking a finite
morphism $\kappa:\check T \to \wt T$, branched over $\wt B$ the
local monodromy groups become unipotent. We consider
\begin{equation}
\begin{picture}(500,500)
\setsqparms[1`1`1`1;400`400] \putsquare(100,100)[\check\cX_m`\wt
X_m`\check T`\wt T;\kappa'`\check f_m`\wt f_m`\kappa]
\end{picture}\qquad .
\end{equation}
The canonical extension  of $R^n \check f'_{m*} \mathbb C_{\check
\cX_m'}\otimes \cO_{\check T'}$ to $\check T$ (\cite{delig}) is a
coherent sheaf. By a theorem of W.~Schmid \cite{schmid}, the
subsheaf $\check f'_{m*}K_{\check \cX'_m/\check T'}$ extends to a
locally free sheaf on $\check T$. Kawamata's theorem \cite{kawa}
states that this locally free extension is equal to $\check f_{m*}
K_{\check \cX_m}$. It is known that also $\wt f_{m*}K_{\wt
\cX_m}$ is locally free: Namely as $\kappa'$ is a proper
holomorphic map of equidimensional complex manifolds, $K_{\wt
\cX_m}\subset\kappa'_*K_{\check\cX}$ is a direct summand, and
hence $\wt f_{m*} K_{\wt \cX_m} \subset \wt f_{m*}
\kappa'_*K_{\check\cX_m} = \kappa_* \check f_{m*}K_{\check\cX_m}$
is a direct summand. Now the latter is locally free, as $\check
f_{m*}K_{\check \cX_m}$ is a locally free $\cO_{\check
T}$-module, and $\kappa$ is a finite proper map of complex
manifolds. We have $K_{\cX_m}= \nu'_*K_{\wt \cX_m}$ on the
manifold $\cX_m$ so that $f_{m*}K_{\cX_m}$ is locally free.

Next, we use W.~Schmid's description of sections of $\check
f_{m*}K_{\check \cX_m}$ around points of the normal crossing
divisor. Let $\Delta^k\simeq U \subset \check T$ be an open subset
such that the complement of the normal crossing divisor is
$U'\simeq \Delta^{*\ell}\times \Delta^{k-\ell}$.

Let $\phi$  be a section of $K_{\check \cX_m}$ over $\check
f_m^{-1}(U)$. Over $U'$ it can be expressed in terms of a basis
$\{s_1,\ldots, s_M\}$ of multivalued (locally constant) sections
of $R^n\check f_{m*}\mathbb C_{\check \cX_m}$ over $U'$. So
$\phi= \sum f_\nu \cdot s_\nu$ for certain multivalued
holomorphic functions on $U'$. According to
\cite[(4.17)]{schmid}, the holomorphicity of $\phi$ in points of
the normal crossing divisor is equivalent to the $f_\nu$ having
at most logarithmic singularities. Next the $L^2$-norm is
computed at points $t\in U'$. (We identify $K_{\check\cX}$ with
$K_{\check\cX/\check T}$).
$$
\|\phi(t)\|^2 = \int_{\check \cX_{m,t}} \phi(t)\wedge
\ol{\phi(t)}= \sum_{\nu\mu}f_\nu(t)\ol{f_{\mu}(t)}\cdot
\int_{\check \cX_{m,t}}s_\nu(t)\wedge\ol{s_{\mu}(t)}
$$
The latter integrals are independent of $t$, because the $s_i$
are locally flat sections.  So
$$
\|\phi\|^2 \leq \sum_{j=1}^\ell c_j (-\log|t_j|)
$$
for some constants $c_j > 0$.

The unipotent reduction preserves such estimates so that a
similar estimate (with different constants) also holds for
sections of $\wt f_{m*}K_{\wt\cX_m}$.

This implies an estimate for sections of $f_{m*}K_{\cX_m}$. We
only note the following rough estimate: Let $W\subset T$ be an
open subset. Then for any $\psi \in (f_{m*}K_{\cX_m})(W)$ we have
$$
\|\psi(x)\|^2 \leq \sum \alpha_j(-\log|\tau_j|)
$$
for certain positive constants $\alpha_j$ and holomorphic
functions $\tau_j$, which vanish on $\cA$. This proves the
following lemma:

\begin{lemma}\label{lelonglemma}
The holomorphic line bundle $\det(f_{m*}K_{\cX_m/T})$ carries a
singular hermitian metric $h$, which is of class $C^\infty$ on
$T\backslash \cA$ such that in local holomorphic coordinates
\begin{equation}\label{logest}
h \leq \sum_j \beta_j (-\log|\tau_j|)
\end{equation}
for certain $\beta_j>0$.
\end{lemma}

\begin{corollary}
For any $x \in T=P\times S$
$$
\sup\{ \gamma \geq 0; -\log h(z) \leq \gamma\cdot \log(\|z-x\|^2)
+ O(1)\}=0.
$$
\end{corollary}

In Proposition~\ref{curvper} we showed that $h$ is a {\it
positive} singular hermitian metric so that the corollary implies
the theorem for $p=1$.

The curvature form $\Theta$ satisfies a Poincar\'e growth
condition on $\Delta^{*\ell}\times \Delta^{k-\ell}$ (cf.\
Lemma~\ref{delta}), i.e.\
$$
\Theta \leq C\cdot \sqrt{-1}\left(\sum_{j=1}^\ell
\frac{dt_j\wedge\ol{dt_j}}{|t_j|^2\log^2|t_j|}+\sum_{j=\ell+1}^k
dt_j\wedge\ol{dt_j}\right).
$$
In particular all powers $\Theta^p$ define closed
$(p,p)$-currents. These estimates hold for the Hodge metrics over
$\check T$, $\wt T$, and since $\wt T \to T$ is a modification of
complex manifolds, also the $\Theta_h^p$ on $T$ are closed
currents.

We show the last statement of Theorem~\ref{lelong}.

Let $x\in P\times S$ be a point and $z_1,\ldots,z_k$ local
coordinates such that $x=0$. Let (locally) $h= {\bf e}^{-u}$,
with $u$ plurisubharmonic, and define $\varphi = \log\|z\|^2$.
For any positive $(p,p)$-current $R$ and small $r>0$ the quantity
$\nu(R,x,r)$ is defined by
$$
\nu(R,x,r)= \frac{1}{r^{2p}}\int_{\|z\|<r} R \wedge
(dd^c\|z\|^2)^{k-p},
$$
and in terms of Demailly's generalized Lelong numbers
$$
\nu(R,x,r)=\nu(R,\varphi,\log r),
$$
where
$$
\nu(R,\varphi,t)= \int_{\varphi(z)<t} R \wedge
(dd^c\varphi)^{k-p}.
$$
In a straightforward way a generalized Jensen formula can be
proved:
\begin{gather*}
\int_{r_0}^{r_1} \nu((dd^c u)^p,\varphi,t)dt = \int_{\varphi=r_1}
u
(dd^c u)^{p-1}\wedge d^c \varphi \wedge (dd^c \varphi)^{k-p}  \\
\hspace{5cm} -\int_{\varphi=r_0} u (dd^c u)^{p-1}\wedge d^c
\varphi \wedge (dd^c \varphi)^{k-p}  \\ \hspace{5cm} -
\int_{r_0<\varphi<r_1} u (dd^cu)^{p-1}\wedge (dd^c\varphi)^{k-p+1}
\end{gather*}
It is known that for any fixed $r_1$
$$
\nu((dd^c u)^p,x)= \lim_{r\to -\infty}\left(-
\int_{r}^{r_1}\nu((dd^c u)^p,\varphi,t)dt /r \right).
$$
So we need only to show that the right hand side of the above
Jensen formula is bounded for fixed $r_1$, independent of $r$.
This follows immediately, because
\begin{itemize}
\item[(i)] $u \geq - c \log(\sum_j \beta_j (-\log|\tau_j|) )$ by
Lemma~\ref{lelonglemma}
\item[(ii)] As a plurisubharmonic function $u$ is (locally)
bounded from above
\item[(iii)] $dd^c u$ satisfies a Poincar\'e growth condition on
$\wt T$.
\end{itemize}

\section{Convergence property of generalized \pw metrics}
We include the definition of generalized \pw metrics, which can
also be part of a conceptual approach to an analytic theory of
moduli spaces, relating canonically polarized manifolds to framed
manifolds and $m$-framed manifolds. However, we will use later
families of $m$-framed manifolds and Galois coverings instead, and
properties of the period map.

In the first place, generalized \pw metrics are intrinsically
defined K"ahler metrics on the base spaces of universal
deformations. Due to functoriality these will be seen to descend
to moduli spaces.

In this section, we will assume that for all $\varepsilon \in
\mathbb Q$, $0 \leq\varepsilon \leq  \varepsilon_0$ the divisor
$$
K_X + (1-\varepsilon) D
$$
is positive. This condition is satisfied for
$\varepsilon_0=1/m_0$ in our basic situation, where $(X, D)$ is
$m_0$-framed and $D$ positive. The methods of
\cite{tsuji1,rkob1,rkob2,tian-yau} yield unique \ke metrics
$\eta_{X, m}$ on the $V$-manifolds $\wt X_m$ (cf.\
Section~\ref{frV}) of  Ricci-curvature $-1$. As in the smooth
compact case we can see that the $V$-\ke metrics define
generalized \pw metric on the moduli space of framed manifolds as
follows:

Let $\cD \hookrightarrow \cX \to S$ define an effective
holomorphic family of framed manifolds $(\cX_s, \cD_s)_{s \in
S}$. Let $(X, D)=(\cX_{s_0}, \cD_{s_0})$. Let $m \geq m_0$, and
let $\wt X_m$ be equipped with the \ke metric $\eta_{X, m}$. For
any $v \in T_{s_0}S$ denote by
$$
A_{m,v}= A_{m,v,\bar\beta}^\alpha \frac{\partial}{\partial
z^\alpha} dz^{\bar\beta} \in \Gamma(X,{\mathcal A}^{V,
1}_X({\mathcal T}^V_X))
$$
the representative of the \ks class of $v$ according to
Remark~\ref{infdef} in $T^1(X, D)$, which is harmonic with respect
to $\eta_{X, m}$.

\begin{definition}
Let $v,w \in T_{s_0}S$, and $A_{m,v}$, $A_{m,w}$ corresponding
harmonic \ks forms. Then the \pw inner product is
$$
<v,w>_{PW} = \int_{X}  <A_{m,v},A_{m,w}> \omega_{X,m}^n.
$$
\end{definition}

The K"ahler property of the induced form $\omega_{PW,m}$ on $S$ can
be shown in the same way as for the case of smooth, canonically
polarized varieties. Also a fiber integral formula holds for the
\pw form, and a line bundle equipped with a Quillen metric can be
constructed, whose curvature form equals $\omega_{PW,m}$ up to a
constant \cite{bgs}.

On the other hand the tangent cohomology $T^1(X, D)$ can be
computed in terms of the complete \ke metric $\omega_{X'}$ on
$X'=X\backslash D$ as the $H^1_{(2)}(X', \cT_{X'})$ the
$L^2$-cohomology group of the sheaf of holomorphic vector fields
$\cT_{X'}$ \cite{sch}. The $L^2$-structure on the tangent
cohomology defines a \pw metric $\omega_{PW, fr}$ on $\cM_{fr}$.

Let $\Omega_{\cX/S}$ be the relative volume form, i.e.\ a
hermitian metric on $\bigwedge^n \cT_{\cX/S}$, induced by all
$\eta_{\cX_s, m}$, and denote by $\eta_{\cX, m}$ the negative of
its curvature form on the total space. Its restrictions to all
fibers are the \ke forms on the fibers. Let $v =
\partial/\partial s \in T_{s}S$ be a tangent vector, and
$$
\frac{\partial}{\partial s} + a^\alpha \frac{\partial}{\partial
z^\alpha}
$$
the horizontal lift with respect to $\eta_{\cX,m}$. Also in the
case of $V$-structures, its exterior derivative restricted to the
fiber $X$
$$
\ol\partial (a^\alpha )= \frac{\partial a^\alpha}{\partial
z^{\bar\beta}} \frac{\partial}{\partial z^\alpha} dz^{\bar\beta}
$$
equals the harmonic \ks form $A_{m,v}$. For a more detailed
discussion of the \pw inner product and \pw forms for singular
base spaces cf.\ also \cite{fs}.

Denote by $\eta_{\cX_s}$ the usual \ke metrics, and by
$\eta_{\cX}$ the negative of its Ricci form on the total space.

Measuring convergence in $C^{k, \alpha}(X')$-spaces with respect
to quasi-coordinates on $X'=X\backslash D$ the $\eta_{X, m}$ tend
to the complete \ke metric $\omega_{X'}$ on $X'$ \cite{tsuji2}.
In a holomorphic family of framed manifolds, this convergence
yields a convergence of the relative volume forms $\Omega_{\cX/S,
m}$ to the relative volume form $\Omega_{\cX'/S}$ of the smooth
\ke metrics in the spaces $C^{k, \alpha}(\cX')$, $\cX'=\cX
\backslash \cD$. Together with the above fact about the
characterization of harmonic \ks forms we see immediately that
the harmonic \ks forms $A_{m, v}$ converge to the harmonic
$L^2$-integrable \ks forms $A_{fr, v}$ on $X'$ with respect to
the complete \ke metrics on $X'$.

\begin{proposition}
For the generalized \pw metrics on moduli spaces of framed
manifolds
$$
\lim_{m \to\infty} \omega_{PW,m} = \omega_{PW,fr}
$$
holds in any $C^k$-topology. Let $m$ be fixed $\cD
\hookrightarrow \cX \to S_{m, fr}$ be a local universal
holomorphic family of $m$-framed manifolds and $\cX_m \to \cX \to
S_{m, fr}$ the induced family of branched coverings with $\cX_{m,
s}$ canonically polarized such that $S_{m, fr}$ embeds into a
base of a universal family of canonically polarized manifolds,
giving rise to $\kappa: S_{m,fr} \hookrightarrow S_c$, where
$S_c$ carries the usual \pw form $\omega_{PW, can}$. Then
$$
\frac{1}{m}\kappa^*(\omega_{PW,can}) = \omega_{PW,m}.
$$
\end{proposition}
We show the second claim: We have the $V$-structures on the fibers
$\cX_s$, and map the usual \ke metrics to \ke $V$-metrics on the
quotients $\cX_{m, s}/\mathbb Z_m$. Any harmonic \ks $V$-form
lifts to a harmonic \ks form on $\cX_{m, s}$. the factor $1/m$ is
due to the integration over $m$ sheets as opposed to the
integration over the $V$-manifold. \qed

\section{Moduli spaces of framed manifolds}\label{modfr}
In the analytic case, a {\it polarization} of a framed manifold
$(X, D)$ is the assignment of a K"ahler class $\lambda_X \in H^2(X,
\mathbb R)$. Polarizations, which are images of integer valued
cohomology classes, coincide with {\it inhomogeneous}
polarizations in the sense of Mumford (cf.\ \cite{mumfo}). (Here,
we can also allow rational coefficients and consider $\mathbb
Q$-divisors.)

The following definition is also sensible for inhomogeneously polarized
framed projective varieties $(X,D,\lambda_X)$ (over $\mathbb C$).

\begin{definition}
\begin{itemize}
\item[(i)] A compact K"ahler manifold $X$ is called uniruled
over a smooth divisor $D$, if there exists a surjective
meromorphic map $\varphi : \mathbb P_1 \times Y \to X$ with the
following properties: The map $\varphi$ does not allow a
meromorphic factorization over $pr_2 : P_1 \times Y \to Y$. The
restriction of $pr_2$ to the proper transform of $D$ under
$\varphi$ is a modification.
\item[(ii)] A polarized framed manifold $(X, D, \lambda)$ is called
non-uniruled, if the K"ahler manifold $D$ is non-uniruled, and if
$X$ is not uniruled over $D$.
\end{itemize}
\end{definition}

In the analytic category, the (coarse) moduli space of
non-uniruled polarized K"ahler manifolds exists.

For non-uniruled, polarized, projective framed manifolds $(X, D,
\lambda_X)$, first the Hilbert polynomials $P(x)$ for $\lambda_X$
on $X$ and $Q(x)$ for $\lambda_X|D$ are of interest. (If the
polarization $\lambda_X$ is represented by $D$, we have
$Q(x)=P(x)-P(x-1)$).

Let $\lambda_X$ be represented by a basic polar divisor and
corresponding ample line bundle $\cL_X$. As usual, Matsusaka's
big theorem (\cite{mats, liemum}) is applied to $(X, \cL_X)$:
There exists an integer $c>0$ only depending on $P(x)$, such that
for all $m \geq c$ the sheaves $m \cdot \cL_X$ are very ample.

\begin{theorem}\label{alsp} There exists an algebraic space $\cM_{fr}$
in the sense of Artin, which is the coarse moduli space of isomorphism
classes of non-uniruled, polarized, framed projective manifolds $(X, D,
\lambda_X)$ with fixed Hilbert polynomials $P(x)$and $Q(x)$.
\end{theorem}
As non-uniruledness is an open and closed condition for polarized
varieties, we can also impose the condition that both $X$, and
$D$ are non-uniruled. Then the assignment
$(X,D,\lambda_X)\mapsto(X,\lambda_X)$ (with Hilbert polynomials
fixed) defines a natural map $\cM_{fr} \to \cM$ of algebraic
spaces, where $\cM$ denotes the moduli space of uniruled polarized
manifolds. If the divisors $D$ are very ample and represent the
polarization $\lambda_X$ (and $X$ is non-uniruled), $D$ may also
be singular giving rise to a moduli space $\wh{\cM}$ equipped
with a natural morphism $\wh\nu : \wh\cM \to \cM$.

\begin{proof} First, $c>0$ as above is taken and $m \geq c$ fixed and
for all polarized varieties $X$ with Hilbert polynomial $P(x)$ a
corresponding projective embedding $X \hookrightarrow \mathbb P_N$
induced by global sections of $m \cdot \cL_X$ considered. As
subvarieties of $\mathbb P_N$ these $X$ have $P(m\cdot x)$ as Hilbert
polynomials. We denote by $Hilb^P_{\mathbb P_N}$ the Hilbert scheme of
all subvarieties with $P(m\cdot x)$ in the sense of Grothendieck
\cite{groth}. The locus ${\cH} \subset Hilb^P_{\mathbb P_N}$ of all
smooth subvarieties is quasi-projective. Let
\begin{equation}\label{hilbsch}
\settriparms[1`1`1;350] \qtriangle[{\cX}`{\cH} \times \mathbb P_N`
{\cH};i`f`pr_1]
\end{equation}
be the universal flat family. Here the fibers ${\cX}_{s}
= f^{-1}(s)$ for $s \in {\cH}$ carry the
polarization $\cO_{{\cX}_{s}}(1)= m \cdot
\cL_{{\cX}_{s}}$

Next we fix the Hilbert polynomial $Q(x)$ with respect to $D$ and
$\cL|D$. Again, by (\cite{groth}), Th\'eor\`eme~3.1 we are looking
at a functor represented by a projective, flat ${\cH}$-scheme
$\wt\nu:\widehat{\cH}\to \cH$ equipped with a universal flat
family ${\cD} \to {\wh\cH}$. The locus $\cH_{fr}$ of smooth
divisors $\widehat\cH \supset {\cH}_{fr} \stackrel{\pi}{\to}
{\cH}$ is a quasi-projective variety. Explicitly, let $P$ be the
dual of $\mathbb P_N$, then $\cH_{fr} \subset \wh\cH=\cH \times
P$ is a Zariski open subspace. We have
\begin{equation}\label{hilbdia}
\begin{picture}(1050,500)
\settriparms[1`1`1;350]
\putqtriangle(100,50)[\cD`{\wh\cX}`\wh\cH;j``\wh{f}]
\setsqparms[1`1`1`1;350`350]
\putsquare(450,50)%
[\phantom{{\cX}}`{\cX}`\phantom{{\cH}}`{\cH};``f`\wt\nu].
\end{picture}
\end{equation}
The graph $\wh \Gamma \subset \wh\cH\times\wh\cH$ of the
equivalence relation identifying embedded manifolds with singular
framings is mapped properly to the graph $\Gamma\subset \cH
\times \cH$, which defines the moduli space $\cM$ of polarized
projective manifolds. By assumption, the natural map $\Gamma \to
\cH$ is proper, so $\wh \Gamma$ also defines a proper equivalence
relation. This ensures the existence of a natural complex
structure on $\wh\cM$. (Observe that this statement can also be
proved in the non-reduced category). Finally $\wh\cM$ carries the
structure of an algebraic space. The construction is compatible
with the restriction to $\cH_{fr}$. If the above equivalence
relations are given by the action of $G=PGL(N+1,\mathbb C)$ on
$\wh\cH$ and $\cH$ resp.\ the moduli spaces $\wh \cM$, $\cM_{fr}$
and $\cM$ are eventually geometric quotients. In the analytic
case the statement of the Matsusaka-Mumford theorem  is also
valid (cf.\ \cite{sch}) for framed polarized manifolds.

Later we will consider compactifications of the algebraic spaces
$\cM$ and $\wh\cM$ by normal crossings divisors with a morphism
$\ol{\wh\cM}\to\ol\cM$. We can assume that it is induced by a flat
morphism $\ol{\wh\cH} \to \ol\cH$ of compactified Hilbert schemes
of the similar type. \end{proof}

The moduli space $\cM$ is induced by a smooth family of the form
(\ref{hilbsch}) with hyperplane section $\cD'\subset \cX$, such
that the very ample divisors $\cD'_s$ represent a fixed multiple
of the polarizations on $\cX_s$. Let $n=\dim \cX_s$ as before.
According to Fujita's theorem \cite{fujita}, the divisors
$K_{\cX_s} + m \cD'_s$ are ample for $m \geq n+2$. We fix $m>\dim
\cX_s + 3 $ and represent $m [\cD'_s]$ by all possible divisors
$\cD_s$. This gives rise to a diagram of the form (\ref{hilbdia}).
We pull back the divisor $\cD'$ to $\wh X$ and obtain a bundle
space $\underline E' \to \wh\cX$. Let $\underline E \to \wh\cX$ be
the bundle associated to $\cD$. Like in Section~\ref{cycov} we
construct a family of cyclic coverings $f_m: \cX_m \to \wh \cH$
and a diagram
\begin{equation}\label{xm}
\settriparms[1`1`1;350] \qtriangle[\cX_m`\wh\cX`\wh \cH;
\mu`f_m`\wh f]
\end{equation}
where the branch locus of $\mu$ is $\cD\subset\wh\cX$. The fibers
$\cX_{m,s}$ are smooth for $s \in \cH_{fr}$.

\begin{proposition}\label{Xm}
The above construction gives rise to a proper and finite morphism
of algebraic spaces $ \kappa :\cM_{fr} \to \cM_c $ from $\cM_{fr}$
to a component of the moduli space of canonically polarized
(smooth) varieties.

Let $(X,D)$ be a fixed framed manifold with branched covering
$X_m \to X$ as above, and let $\wt R$ and $R$ resp.\ denote base
spaces of universal deformations. Then there exists a closed
holomorphic embedding
$
\wt \kappa : \wt R \to R
$
which induces the map $\kappa$ in a neighborhood of the
corresponding moduli point, where it is of the form
$
\wt R/{\rm Aut}(X,D) \to  R/ {\rm Aut}(X_m).
$
\end{proposition}

\begin{proof} We show the properness of $\kappa$. Let $(X_j, D_j)$, $j \in \mathbb
N$ be framed manifolds such that the images under $\kappa$ of the
isomorphism classes are represented by canonically polarized
manifolds $X_{m, j}$. Assume that the isomorphism classes of
$X_{m, j}$ converge to the class of some canonically polarized
manifold $Z$. We consider a universal deformation ${\mathcal Z}
\to S$ of $Z$ and can assume that the $X_{m, j}$ are fibers. The
space of automorphisms ${\mathcal Aut}_S(\mathcal Z) \to S $ is
mapped properly to the base $S$ with reduced, $0$-dimensional
fibers. Hence, for a subsequence of the $j$, the group of deck
transformations of the $X_{m, j}$ converges to a subgroup of
${\rm Aut(Z)}$. Its fixed point set is a smooth divisor, and a
subsequence of the $(X_j, D_j)$ converges to a framed manifold
$(X, D)$. The bundle $K_X + [D]$ is positive, since $K_Z$ is
positive.

Let $(X,D)$ be a general framed manifold, and $(\wt R, \wt r_0 )$
the base of a universal deformation with fibers $(\cX_r, \cD_r)$.
We have $[D]= m \cdot [D']$ by assumption. As for all $r$ the
$[\cD_r]$ are divisible by $m$, the given family of polarized
manifolds gives rise to a family $\cX_m \to \wt R$ with Galois
covering $\cX_m \to \cX$. This family is induced by a universal
deformation $\cY \to R$ of $X_m$. The base change $\wt R \to R$
must be an embedding by Remark~\ref{infdef}.

Let $\mathbb Z_m \subset {\rm Aut}(X_m)$ be the group of deck
transformations. It acts on $R$ leaving the subspace $\wt R
\subset R$ pointwise fixed, since the group action can be lifted
to all of its fibers. Furthermore, there is a natural morphism
${\rm Aut(X, D)} \to {\rm Aut}(X_m)/\mathbb Z_m$ so that we have a
holomorphic map $\wt R/{\rm Aut(X, D)} \to \wt R / ({\rm
Aut}(X_m)/\mathbb Z_m) \subset R / {\rm Aut}(X_m)$, which is
given by the action of a finite group.
\end{proof}

\section{Fiber integrals and determinant line bundles for morphisms}
We will use the method of {\it generalized determinant line bundles}.
Let $F: Z \to S$ be a proper, holomorphic map of complex spaces and
$\cL$ a coherent $\cO_Z$-module.

The direct image $R^\bullet F_* {\mathcal L}$ of ${\mathcal L}$ under
the proper map $F$ in the derived category can be locally represented
by a sequence ${\mathcal F}^\bullet$ of finite, free ${\mathcal
O}_S$-modules, which is bounded to the right. If the morphism is flat,
the sequence can be chosen as bounded, and the tensor product of the
determinant sheaves of the ${\mathcal F}^i$ with alternating exponents
$\pm 1$ is by definition the determinant line bundle $\lambda =
\det({\mathcal L})$, and the latter is globally well-defined.

Let $\cL=\cO_Z(L)$ be a holomorphic line bundle equipped with a
hermitian metric of class $C^\infty$. According to Bismut, Gillet
and Soul\'e, \cite{bgs}, under the assumption of $F$ being a
smooth K"ahler morphism of complex manifolds (or reduced complex
spaces \cite{fs}), the Chern form of the Quillen metric $h^Q$ on
$\det({\mathcal L})$ is equal to the  component of degree two of a
fiber integral:
\begin{equation}\label{BGS}
c_1({\lambda}, h^Q) = -\left[\int_{Z/S} {\rm td}(Z/S){\rm
ch}({\mathcal L})\right]_{(2)},
\end{equation}
where $\rm td$ and $\rm ch$ resp.\ define the Todd and Chern character
resp. (This holds also, when $\mathcal L$ is replaced by a hermitian
vector bundle.)

By functoriality and universal properties, this equation extends
to ${\mathcal L}$ replaced by an element of the Grothendieck
group, i.e.\ a virtual holomorphic vector bundle. If $n$ denotes
the fiber dimension, the virtual bundle $({\mathcal L}-{\mathcal
L}^{-1})^{n+1}$ has rank zero, and the lowest term in ${\rm
ch}(({\mathcal L}-{\mathcal L}^{-1})^{n+1})$ is
$2^{n+1}c_1({\mathcal L})$ so that the only contribution of the
Todd character in (\ref{BGS}) is equal to $1$. Hence the Chern
form of $\det(({\mathcal L}-{\mathcal L}^{-1})^{n+1})$ equals
\begin{equation}\label{fint1}
- 2^{n+1} \int_{Z/S} c_1({\mathcal L},h)^{n+1}.
\end{equation}
Now we return to the situation of moduli spaces like in
Section~\ref{modfr}. The Hilbert scheme $\cH_{fr}$ carries the
determinant line bundle $\lambda_{fr}$ with singular hermitian
metric $h_{fr}$ according to Proposition~\ref{curvper} and
Lemma~\ref{lelong}. It is important that the line bundle
$\lambda_{fr}$ on $\cH_{fr}$ was extended to the line bundle $\wh
\lambda $ on $\wh\cH$. Let $\pi: \cH_r \to \cH$ be a
desingularization with fiber product $\nu_r:\wh\cH_r \to \cH_r$
and pull-back $\wh \lambda_r$ of $\wh\lambda$. Since $\nu_r$ is a
smooth map, we can apply the above methods. We consider the
determinant bundle $\det((\wh\lambda_r -
\wh\lambda_r^{-1})^{N+1})$.

We now apply these methods to singular hermitian metrics on
singular spaces (cf.\ Section~\ref{singherm}), and
$(1,1)$-currents.

So far we are given a smooth holomorphic map $\wh\nu : \wh\cH \to \cH$
and a holomorphic line bundle on $\wh\lambda$ on $\wh\cH$, whose
restriction $\lambda_{fr}$ to $\cH_{fr}$ carries the $C^\infty$
hermitian metric $h_{fr}$ with curvature form $\Theta_{fr}$.

We use the above arguments to extend the determinant line bundle
$\det((\wh\lambda -\wh\lambda^{-1})^{N+1})$ as a coherent sheaf
from $\cH$ to $\overline{\cH}$. We denote by $\wh\Theta$ the
curvature current of $\wh\lambda$. Let $\ell=\dim \cH$. In order
to define a fiber integral
$$
\int_{\wh\cH/\cH} \wh\Theta^{N+1},
$$
for any $(\ell-1,\ell-1)$-form $\varphi$ of class $C^\infty$ with
compact support, we set
$$
\Theta^Q(\varphi) =\int_{\wh\cH/\cH} \Theta_{fr}^{N+1} \wedge
\wt\nu^*\varphi,
$$
with $\Theta_{fr}= \wh\Theta|\cH_{fr}$.

At this point, we may blow up $\wh\cH$ with exceptional set in
$\wh\cH\backslash\cH_{fr}$ and realize $\cH_{fr}$ as a complement
of a divisor with only normal crossings singularities so that the
assumptions of Lemma~\ref{delta} are satisfied. The upper
Poincar\'e growth estimate for $\Theta_{fr}$ implies that the
above integral is finite, and it vanishes, if $\varphi$ is
$d$-exact. So $\Theta^Q$ is well-defined as a $d$-closed
$(1,1)$-current. Also Lemma~\ref{delta} implies that $\Theta^Q$
is positive (in the sense of currents).

\begin{proposition}\label{lelofib}
At all points $\cH$ the Lelong numbers of $\Theta^Q$ vanish.
\end{proposition}
The above statement also holds after descending to the moduli
space at points of the boundary, as we can always achieve the
situation of Section~\ref{singherm} after blowing up the boundary.
\begin{proof}
The proof follows immediately from Theorem~\ref{lelong}.
\end{proof}
\begin{lemma}\label{chern}
The current $(1/2\pi)\Theta^Q$ on $\cH$ represents the Chern-class of
the bundles $\det((\wh\lambda-\wh\lambda^{-1})^{N+1})$ on $\cH$.
\end{lemma}

\begin{proof} We use an auxiliary $C^\infty$ hermitian metric $h_a$ on
$\wh\lambda$ with curvature form $\Theta_a$. Then the fiber
integral $\int_{\wh\cH/\cH} \Theta_a^{N+1}$ exists and represents
up to a numerical constant the Chern class $c_1(\det((\wh
\lambda-\wh \lambda^{-1})^{N+1}))$ on $\cH$. On $\cH_{fr}$ the
difference $\Theta_{fr}-\Theta_a$ is (globally) of the form
$\sqrt{-1}\partial\ol\partial u$. Now
$$
\Theta_{fr}^{N+1}= \sqrt{-1}\partial\ol\partial u\wedge \Omega +
\Theta_a^{N+1},
$$
where $\Omega= \sum_{j=0}^N \Theta_{fr}^j
\wedge \Theta_a^{N-j}$.

Basic properties of the $L^2$-Dolbeault-complex on
$\Delta^{*k}\times \Delta^l$ (cf.\ \cite{zuck}) show that $u$, and
$\ol\partial u$ can be chosen as locally $L^2$-integrable (with
respect to metrics with Poincar\'e growth condition). So
$\int_{\cH_{fr}/\cH} \sqrt{-1}\, \ol\partial u \wedge \Omega $
defines actually a current. We claim that in the sense of currents
\begin{equation}\label{Omega}
\int_{\cH_{fr}/\cH} \sqrt{-1}\, \partial \ol\partial u \wedge
\Omega = - d   \int_{\cH_{fr}/\cH} \sqrt{-1}\, \ol\partial u
\wedge \Omega
\end{equation}
holds.

In fact, the right hand side applied to a $C^\infty$-form with compact
support equals
$$
-d\int_{\cH_{fr}} \ol\partial u \wedge \Omega \wedge \wt\nu^*
d\varphi = \int_{\cH_{fr}} \partial \ol\partial u \wedge \Omega
\wedge \wt\nu^* \varphi = \left(\int_{\cH_{fr}} \partial
\ol\partial u \wedge \Omega \right) (\varphi).
$$
\end{proof}
\begin{corollary}There exists a singular hermitian metric $h^Q$ for
$\det((\wh\lambda-\wh\lambda^{-1})^{N+1})$ on $\cH$, whose
curvature is positive in the sense of currents.
\end{corollary}

\begin{remark} Furthermore, it follows from the construction that for
any subspace of $\cH$, in particular for any curve in $\cH$, the
restriction of $h^Q$ and $ \Theta^Q$ resp.\ exist as singular metric and
$d$-closed current resp. If $C \subset \cH$ is a local analytic curve
through a point $p$, representing a direction, where $\cX \to \cH$ is
effective, the current is strictly positive in this direction.
\end{remark}

The latter fact follows immediately, because the form $\Theta_{fr}$ is
strictly positive on the preimage of $C$ in $\cH_{fr}$.

After blowing up the boundary  $\lambda^Q$ possesses a line bundle
extension $\overline{\lambda^Q}$ on $\overline{\cH}$.  The result
of this paragraph concerning Hilbert schemes is so far:

\begin{theorem} The compactified Hilbert scheme $\overline{\cH} \supset
\cH$ carries a line bundle $\overline{\lambda}^Q$ with a singular
hermitian metric $\overline{h^Q}$ whose curvature
$\overline{\Theta^Q}$ is positive. The Lelong numbers vanish
everywhere, and $\Theta^Q$ is strictly positive in effective
directions of the family $\cX \to \cH$. Moreover, on $\cH$ the
construction is functorial with respect to base changes of
families concerning the line bundle and its curvature.
\end{theorem}

In a final step we descend to the moduli space $\overline\cM$.

The automorphism groups of the polarized manifolds act on local
universal deformation spaces in a finite way (with uniformly
bounded orders). By functoriality, a certain power
$(\lambda^Q)^\mu$ descends from $\cH$ to some $\lambda_\cM$ on
$\cM$. On $\overline\cM$ the line bundle $\overline{\lambda^Q}$
gives rise to a coherent sheaf. As $\mu \cdot \Theta^Q$ is
invariant under the action of the projective linear group on
$\cH$, it descends to a current $\Theta_\cM$ on $\cM$.  We look
at the natural map $u: \cH \to \cM$ extended to $\ol
u:\overline{\cH} \to {\overline\cM}$. The current $\Theta_\cM$
will now be extended to $\ol\cM$: Let $\varphi$ be a $C^\infty$
differential form of degree $(\dim \cM -1, \dim \cM -1)$ with
compact support. We take a closed subvariety $S \subset \ol\cH$,
so that the map $S \to \ol\cM$ is generically finite, and
dominant. The following definition is independent of the choice
of $S$:
\begin{equation*}
\begin{split}
\Theta_{\overline \cM}(\varphi)= \int_\cM \Theta_\cM \wedge
\varphi  := \frac{1}{\alpha}\int_S \mu \cdot \Theta^Q \wedge \ol
u^*(\varphi),
\end{split}
\end{equation*}
where $\alpha$ denotes the generic degree of the map $u|S : S \to
\ol\cM$. With $\varphi=d\psi$ we see the closedness of the
current. Again, we have a positive $d$-closed current
$\Theta_{\overline \cM}$. It realizes the Chern class of
$\lambda_\cM$ on $\cM$, which is the restriction of a coherent
sheaf on $\overline{\cM}$. Again, after blowing up the boundary
it possess a line bundle extension $\lambda_{\overline\cM}$ with
a corresponding singular hermitian metric $h_{\overline \cM}$.

\begin{theorem}\label{main1} The moduli space $\cM$ possesses a
compactification $\overline \cM$ as an algebraic space and a
holomorphic line bundle $\overline \lambda $ with a singular
hermitian metric $h$ of positive curvature form $\Theta_h$ such
that
\begin{enumerate}
\item[(i)] for all $p\in \ol\cM$ and any holomorphic curve $C \subset
\ol\cM$ through $p$ with $C\cap \cM \neq \emptyset$ the (positive,
$d$-closed) current $\Theta_h|C$ is well-defined, and the Lelong
number $\nu(\Theta_h|C,p)$ vanishes,
\item[(ii)] for any smooth locally closed subspace $Z \subset \cM$
the current $\Theta_h|_Z$ is well-defined, and $\Theta_h|_Z \geq
\eta_Z$ in the sense of currents, where $\eta_Z$ denotes some
$C^\infty$ hermitian form on $Z$.
\end{enumerate}
\end{theorem}

\section{$L^2$-methods}\label{L2}
In this section, we gather some results based upon H"ormander's
techniques (cf.\ also the result by Ohsawa and Takegoshi
\cite{oh-ta}).

Let $(Y, \omega_Y)$ be a complete K"ahler manifold, and $(L,h)$ be
a hermitian line bundle on $Y$. We write
$$
\omega_Y =
{\frac{\sqrt{-1}}{2}}g_{\alpha\bar\beta}dz^\alpha\wedge
dz^{\bar\beta},
$$
and use the semi-colon notation for covariant derivatives with
respect to the metric tensor. Moreover the components of the
connection form of the line bundle are
$$
\theta_\alpha=\frac{h_{;\alpha}}{h},
$$
and we use $\Theta_{\alpha\ol\beta}=\theta_{\alpha;\ol\beta}$ for
the coefficients of the curvature tensor. We use $\nabla_\alpha$
to denote covariant derivatives of $L$-valued tensors, and
$\|..\|$, $\|..(p)\|$ resp.\ denote norms and pointwise norms
resp. Let $\varphi=\varphi_{\ol\beta}\;dz^{\ol \beta}$ be any
$L$-valued $(0,1)$-form of class $C^\infty$. Then
$$
\overline\partial^*\varphi = - g^{\bar\beta\alpha} \nabla_{\alpha}
\varphi_{\bar\beta} = - g^{\bar\beta\alpha}(\varphi_{\bar\beta;
\alpha} + \varphi_{\bar\beta} \theta_\alpha)
$$
is the formal adjoint of the $\ol\partial$-operator.

The {\em rough Laplacian} is defined by
$$
\Delta\hspace{-10pt}{-} \varphi
= -g^{\bar\delta\gamma} \nabla_\gamma \nabla_{\bar\delta} \varphi_{\bar\beta}
 dz^{\bar\beta},
$$
and the Bochner-Kodaira-Weitzenboeck formula for this situation
reads
$$
\openbox \varphi= (\ol\partial\ol\partial^* +
\ol\partial^*\ol\partial )\varphi = \Delta\hspace{-10pt}{-}
\varphi + g^{\bar\delta\alpha}\varphi_{\ol\delta}
(R_{\alpha\bar\beta} + \Theta_{\alpha\bar\beta}) dz^{\bar\beta} ,
$$
where $R_{\alpha\bar\beta}$ denotes the Ricci-tensor of
$\omega_X$. (The contribution $R_{\alpha\bar\beta}$ cancels out,
if we replace $L$ by $L + K_Y$). The formula implies
\begin{equation}\label{lowest}
\|\ol\partial \varphi \|^2 + \|\ol\partial^* \varphi \|^2 \geq
\int_Y \varphi_{\bar\beta} (R_{\gamma\bar\delta} +
\Theta_{\gamma\bar\delta}) \overline{\varphi}_\alpha
g^{\bar\beta\gamma} g^{\bar\delta\alpha}\; h \; g \; dV_\omega
\end{equation}
for all $C^\infty$-forms with compact support. According to
Andreotti and Vesentini \cite{and-ves} the estimate
(\ref{lowest}) holds (use cut-off functions) for all square
integrable forms $\varphi$, for which $\ol \partial \varphi$ and
$\ol\partial^*\varphi$, taken in the distributional sense, are
square integrable.  Let $H_1$ and $H_2$ resp.\ be the Hilbert
spaces of square integrable $L$-valued $(n,0)$- and $(n,1)$-forms
resp. Then the exterior derivative $\ol \partial$ is a densely
defined closed operator $T: H_1 \to H_2$ whose adjoint $T^*$ is
given by $\ol\partial^*$ (cf.\ \cite{and-ves}).

\begin{proposition}\label{hoer}
Let $(Y,\omega_Y)$ be a K"ahler manifold, which possesses also a
complete K"ahler metric, and let $(L,h)$ be a holomorphic line
bundle, with a singular hermitian metric, whose Lelong numbers
vanish everywhere. Suppose that
$$
\Theta_h \geq c(p) \cdot \omega_Y
$$
for some continuous, everywhere positive function $c(p)$ on $Y$.
Then for any $L$-valued $(n,1)$-current $v$ with $\ol \partial
v=0$, and
$$
\int_Y \frac{1}{c(p)} \|v(p)\|^2 dV_\omega < \infty
$$
there exists an $L$-valued current $u$ with $\ol\partial u= v$
and
$$
\int_Y \|u(p)\|^2 dV_\omega \leq \int_Y \frac{1}{c(p)} \|v(p)\|^2
dV_\omega.
$$
\end{proposition}

\begin{proof} We assume first that $h$ is of class $C^\infty$ and
that $\omega_Y$ is complete. We follow the argument of H"ormander
and Demailly. The closed subspace $F\subset H_2$ of all
$\ol\partial$ closed forms contains the range of $T$, and $T^*$
vanishes on the orthogonal complement of $F$ so that we can
consider $T$ as an operator from $H_1$ to $F$, and $T^*$ as an
operator from $F$ to $H_1$. Now (\ref{lowest}) implies for all
$g\in F$, contained in the domain of $T^*$ with
$\int(1/c(p))\|v(p)\|^2  <\infty$ that
$$
\|T^*(g)\|^2 \geq \int_Y c(p) \|g(p)\|^2 dV_\omega.
$$
For any $g,v \in F$, with $g$ in the domain of $T^*$ and
$\int(1/c(p))\|v(p)\|^2  <\infty$, we have
$$
|(g,v)|^2 \leq \int_Y \frac{1}{c(p)} \|v(p)\|^2 dV_\omega\cdot
\int_Y c(p) \|g(p)\|^2 dV_\omega,
$$
hence
$$
|(g,v)|\leq \left( \int_Y \frac{1}{c(p)} \|v(p)\|^2
dV_\omega\right)^{1/2}\cdot \|T^*(g)\|.
$$
For any such $v$ there is a continuous linear functional on the
range of $T^*$ sending $T^*g$ to $(g,v)$. The Hahn-Banach theorem
implies the existence of some $u \in H_1$ such that
$(T^*g,u)=(g,v)$ for all $g$ in the domain of $T^*$, i.e.\ $v=Tu$.
Moreover $\|u\|^2 \leq\int_Y (1/c(p)) \|v(p)\|^2 dV_\omega$.

This result was extended by Demailly to arbitrary K"ahler metrics
in \cite{demai1}. The generalization to singular hermitian metrics
is due to Nadel \cite{nadel}.
\end{proof}

\section{Multiplier ideal sheaves}\label{mulshea}
Let $(L, h)$ be a singular hermitian line bundle on a complex
manifold $M$. The sheaf $\cL^2(L, h)$ of square-integrable
sections with respect to $h$ is defined by
$$
\cL^2(L,h)(U)= \{ \sigma \in \Gamma(U, \cO_M(L));
h(\sigma,\sigma) \in L^1_{loc}(U)\},
$$
for open subsets $U\subset M$. There exists an ideal sheaf
$\cI(h)$, called {\em multiplier ideal sheaf} such that
$$
\cL^2(L,h)(U)= (\cO_M(L)\otimes \cI(h))(U)
$$
holds. If we write $h={\bf e}^{-\varphi}\cdot h_0$, where $h_0$
is a hermitian metric of class $C^\infty$, and $\varphi \in
L^1_{loc}(M)$ is the weight function, we see that
$$
\cI(h)=\cL^2(\cO_M,{\bf e}^{-\varphi})
$$
holds. We also use the notation $\cI(\varphi)$ for this sheaf.

For any modification $\pi: \wt M \to M$ of complex manifolds, and
any plurisubharmonic function $\chi$ the following identity of
multiplier ideal sheaves is known (cf.\ \cite[Prop.\
5.8]{demai4}):
\begin{equation}\label{comp_mult}
\pi_*(\cO_{\wt M}(K_{\wt M})\otimes \cI(\chi \circ \pi) )=
\cO_{M}(K_{M})\otimes \cI(\chi).
\end{equation}

\begin{definition}\label{analsing}
A plurisubharmonic function $\varphi$ on a complex manifold is
said to have analytic singularities, if locally
$$
\varphi= \alpha \log(\sum_1^k |f_i|^2) + \varphi_0,
$$
where the $f_i$ denote holomorphic functions, $\varphi_0$ is a
$C^{\infty}$-function, and $\alpha \in \mathbb R_+$.

If $\sigma_i$ are global sections of a line bundle $L$,
$$
h^\alpha= \frac{{\bf e}^{-\varphi_0}}{(\sum|\sigma_i|^2)^\alpha}
$$
defines a singular hermitian metric of positive curvature. In the
above sense it will be called a metric with {\it analytic
singularities} or {\it algebraic singularities} resp. (In the
latter case $\alpha \in \mathbb Q_+$ is also required.)
\end{definition}

In the above situation the holomorphic functions $f_i$ define
some ideal $\cJ \subset \cO_M$. We blow up $M$ along the ideal
$\cJ$, desingularize in a way such that the exceptional set of
the blow-up becomes a normal crossing divisor $D=\sum D_I$. We
call the resulting modification $\pi : \wt M \to M$. Now
\begin{equation}\label{canpul}
K_{\wt M} = \pi^* K_M + R,
\end{equation}
where $R= \sum \rho_j D_j$, $\rho_i \in \mathbb N$ is the
exceptional divisor of $\pi$ on $\wt M$.

The pull-back of $\sum|f_i|^2$ to $\wt M$ is of the same type,
with zero-set equal to $D$, so it is of the form $\prod
|\tau_i|^{2\beta_i}(1+ \sum |\wt f_j|^2)$, so
\begin{equation}\label{pulba}
\varphi\circ\pi = \sum \beta_i \log |\tau_i|^{2} + \wt\varphi_0,
\end{equation}
where $\{\tau_{i}\}$ are defining functions of $\{ D_{i}\}$,
$\beta_i \in \mathbb R_{\geq 0}$, and $\wt\varphi_0$ is some
$C^\infty$-function. In this case the multiplier ideal sheaf can
be computed explicitly as
\begin{equation}\label{muliba}
\cI(\varphi \circ \pi)= \cO_{\wt M}(- \sum \lfloor\beta _i
\rfloor D_i ),
\end{equation}
where $\lfloor
\beta _i \rfloor$ is the Gaussian bracket. Together with
(\ref{comp_mult}) this implies
\begin{equation}\label{multid}
\cI(\varphi) = \pi_*\cO_{\wt M} \left(\sum (\rho_i - \lfloor
\beta _i \rfloor) D_i\right).
\end{equation}
In particular, $\cI(\varphi\circ\pi)$ is locally free.

\begin{proposition}\label{id_mult}
 Let $\Delta^n \subset \mathbb C^n$ be
a polydisk, $\varphi$ a plurisubharmonic function with analytic
singularities on $\Delta^n$, and $\psi$ a plurisubharmonic
function such that $\sqrt{-1} \partial \ol \partial \psi$ is
absolutely continuous on any local holomorphic curve $C \subset
\Delta^n$ with $\varphi|C \not\equiv - \infty$.

Then, after replacing $\Delta^n$ by any smaller, relatively
compact polydisk, there exists real numbers $\gamma$ arbitrarily
close to $1$ such that
$$
\cI(\gamma\cdot\varphi ) = \cI(\gamma \cdot \varphi +\psi )
$$
holds.
\end{proposition}
\begin{proof}
In the sequel, we always allow $\Delta^n$ to be replaced by a
slightly smaller polydisk. We first apply the above modification
to $M=\Delta^n$ with respect to $\varphi$. Then we perform a
further sequence of blow-ups and get a modification $\pi: \wt
\Delta \to \Delta^n$ so that also $\cJ=\cI((\varphi +\psi)\circ
\pi )$ is locally free, and such that with $\wt M = \wt \Delta$
the exceptional divisor is of the above form $D=\sum D_i$. We
still have (\ref{pulba},\ref{muliba}) for $\varphi$.

For any point $x \in \wt\Delta \backslash D$ the function
$\varphi \circ \pi$ is of class $C^\infty$, and $\psi \circ \pi$
is absolutely continuous, when restricted to curves through $x$.
Hence the Lelong number $\nu((\varphi + \psi)\circ \pi,x)$
vanishes. By \cite{bom,sko} we have $\cJ_x = \cO_{\wt \Delta,x}$.
So $\cJ$ is locally free with $V(\cJ)\subset D$. Hence $\cJ =
\cO_{\wt \Delta}(-\sum \beta^\prime_i D_i)$ for some nonnegative
integers $\beta_i'$.

Next, we use (\ref{canpul}) as above and get
\begin{equation}
\cI(\alpha \cdot \varphi) = \pi_*(\cO_{\wt \Delta}(\sum(\rho_i -
\lfloor \alpha \beta_i \rfloor) D_i)
\end{equation}
for all $\alpha >0$.

We chose $\alpha$ so that $\alpha\beta_i \not\in \mathbb Z$ for
all $\beta_i\neq 0$. Next, we compute Lelong numbers. Let $x \in
\wt \Delta$ and $\wt C \subset \wt \Delta$ a local analytic curve
through $x$. If $\pi(\wt C)$ is a point, at which $ \psi $ is
different from $-\infty $, the Lelong number of $\psi \circ \pi$
vanishes. If $\pi(\wt C)$ is a curve $C$, the assumption that
$\psi|C$ is absolutely continuous implies that $\nu(\psi\circ
\pi|\wt C, x)=0$. Because of the additivity of Lelong numbers
$\nu(\pi^*(\alpha\varphi + \psi)|\wt C, x) =
\nu(\pi^*(\alpha\varphi)|\wt C, x)$.

So far $\nu(\pi^*(\alpha \varphi+\psi), x) =
\nu(\pi^*(\alpha\varphi), x)$ holds on $\wt\Delta$. For any point
$x \in D_i \backslash\sum_{j\neq i} D_j$ this Lelong number is
equal to $\nu_i:=\alpha \beta_i \not\in \mathbb Z$.

The latter fact allows us to compute the multiplier ideal sheaf
from the Lelong number: As it is locally free and the space is
smooth, it is sufficient to compute it for points on the regular
part of  the normal crossings divisor $D$. Let $D_i$ be the zero
set of a coordinate function $\tau_i$. Then
$$
0\leq \nu( |\tau_i|^{2\lfloor\nu_i\rfloor} {\bf
e}^{-(\alpha\varphi+\psi)\circ \pi },x) <1
$$
at some $x\in D_i\backslash \bigcup_{j\neq i}D_i$. It follows
from \cite{bom,sko} that
$$
\cI(|\tau_i|^{2\lfloor\nu_i\rfloor} {\bf
e}^{-(\alpha\varphi+\psi)\circ \pi} )_x = \cO_{\wt \Delta,x}
$$
i.e.\
$\tau_i^{\lfloor\nu_i\rfloor}\in\cI((\alpha\varphi+\psi)\circ \pi
)_x$. We need to see that no lower power $\tau_i^k$ is contained
in this multiplier ideal sheaf.

From the Lelong number of $h$, we get the known lower estimate
$$
h\geq \frac{C}{\|z-x\|^{2\nu_i}}.
$$
We use this estimate on a local analytic curve $C_x$, which
intersects $D_i$ in $x$ transversally. So $\int_{C_x}
h|\tau_i|^{2k} dV_{C_x} = \infty$. The same argument is used for
points on $D_i$ near $x$. By Fubini's theorem $\tau_i^k$ is not
in the multiplier ideal sheaf. Now equation (\ref{comp_mult})
implies the claim.
\end{proof}
\begin{remark}\label{perturb}
The above proposition is still valid for the wider class of those
\psh\ functions, which differ from a \psh\ function with analytic
singularities, by a function, which is bounded by $c \cdot
\log(-\log \delta(x)))$, where $c>0$ is a constant, and $\delta$
is the distance of $x$ from the singular set.
\end{remark}

\section{A criterion for quasi-projectivity}\label{crit-q}
Let $X$ be a not necessarily reduced algebraic space with
compactification $\ol{X}$ in the sense of algebraic spaces, and
let $L$ be a holomorphic line bundle on $\ol{X}$ with a {\it
positive} singular hermitian metric $h$ in the sense of
Section~\ref{singherm}.
\begin{condition}[\bf P]
We say that the positivity condition (P) holds, if
\begin{enumerate}
\item[(i)] for all $p\in \ol X$ and any holomorphic curve $C \subset
\ol X$ through $p$ with $C\cap X \neq \emptyset$ the (positive,
$d$-closed) current $\Theta_h|C$ is well-defined, and the Lelong
number $\nu(\Theta_h|C,p)$ vanishes,
\item[(ii)] for any smooth locally closed subspace $Z \subset X$
the current $\Theta_h|_Z$ is well-defined, and $\Theta_h|_Z \geq
\gamma_Z$ in the sense of currents, where $\gamma_Z$ denotes some
$C^\infty$ hermitian form on $Z$.
\end{enumerate}
\end{condition}

Now we state the criterion.

\begin{theorem}\label{criterion}
Let $X$ be an irreducible, not necessarily reduced algebraic space
with a compactification $\ol X$. Let $L$ be a holomorphic line
bundle on $\ol X$. The map
$$
\Phi_{|m L|} : \ol X \rightharpoonup \mathbb P_N(m),
$$
where $N(m)=\dim |mL|$, defines an embedding of $X$ for
sufficiently large $m$, if it satisfies condition (P).
\end{theorem}

We will prove the theorem by induction over $n= \dim X$. The case
$n=1$ is obvious: Let $X$ be an algebraic curve. If $\ol X$ is smooth,
the assumption implies that ${\rm deg}(L) >0$. Let $\ol X$ be a
singular curve and $\pi :\wt X \to \ol X$ be the normalization. Then
${\rm deg}(\pi^* L) >0$ from the assumption so that $L^{\otimes \ell}$
defines an embedding of $\ol X$ into a projective space.

\section{Bigness of $L$}
Let $X$ be a reduced, compact complex space of dimension $n$, and
$\cL =\cO_X(L)\in Coh(X)$ an invertible sheaf.
\begin{definition}\label{defbig}
The sheaf $\cL$ is called  big, if
$$
\limsup_{m \to \infty} \frac{1}{m^n} h^0(X, \cL^{\otimes m})>0.
$$
\end{definition}
In the sequel we denote by $\nu : Y \to X$ The normalization of
the (not necessarily locally irreducible) space $X$, and by $\rho
: Z \to Y$ a modification such that $Z$ is smooth. If $X$ is a
Moishezon space, we assume also that $Z$ is projective. Let $\pi=
\rho \circ \nu$.

\begin{proposition}\label{equivbig}
The following are equivalent:
\begin{itemize}
    \item[(i)] $\cL$ is big
    \item[(ii)] $\nu^*\cL$ is big
    \item[(iii)] $\pi^*\cL$ is big
\end{itemize}
\end{proposition}
\begin{proof}
We show that (ii) implies (i): Consider the exact sequence of
$\cO_X$-modules
$$
0 \to \cO_X \to \nu_* \cO_Y \to \cC \to 0,
$$
where ${\rm supp}(\cC) \subsetneq X$ is nowhere dense, and
$$
0 \to \cL^{\otimes m} \to \nu_*\nu^* \cL^{\otimes m} \to
\cC\otimes {\cL^{\otimes m}} \to 0.
$$
The claim follows, because $h^0({\rm supp}(\cC), \cL^{\otimes
m}\otimes \cC)= O(m^{n-1})$, and $h^0(Y,\nu^*\cL^{\otimes m})\sim
m^n$ .

The other implications are obvious.
\end{proof}

For any $m>0$ with $h(X, \cL^{\otimes m})>0$ we denote by
$\Phi_{\cL^{\otimes m}} : X \to \mathbb P_N$, $N=N(m)$ the
meromorphic map induced by global sections.

\begin{proposition}\label{ZOemb}
Let $X$ be a (reduced) Moishezon space then the following are
equivalent:
\begin{itemize}
    \item[(i)] $\cL$ is big
    \item[(ii)] $\Phi_{\cL^{\otimes m}}: X \to \mathbb P_N$
    embeds some Zariski open subset of $X$.
    \item[(iii)] $\dim \Phi_{\cL^{\otimes m}}(X)= \dim X$
\end{itemize}
\end{proposition}
\begin{proof}
We need to show that (i) implies (ii), the remaining implications
are clear.

We consider as above the normalization and desingularization maps
with $Z$ projective. By Proposition~\ref{equivbig}, $\pi^*\cL$ is
big on $Z$. Let $\cA$ be a very ample invertible sheaf. By
Kodaira's lemma (cf.\ \cite[App.]{k-o}), for some $m>0$ the sheaf
$\pi^*\cL^{\otimes m}\otimes \cA^{-1}$ possesses a non-zero
section with zero divisor $E$ so that the sections of
$\pi^*\cL^{\otimes m}$ yield an embedding of $Z\backslash E$ into
some $\mathbb P_N$. As $H^0(Z,\pi^*\cL^{\otimes
m})=H^0(Y,\nu^*\cL^{\otimes m})$, the invertible sheaf
$\nu^*\cL^{\otimes m}$ gives rise to an embedding of some Zariski
open subset of $Y$ into $\mathbb P_N$. Consider
$$
0 \to \cO_X \to \nu_*\cO_Y \to \cC \to 0.
$$
Let $\cI \subset \cO_X$ be the annihilator of $\cC$. The zero set
$V(\cI)\subset X$, consisting of all non-normal points of $X$, is
nowhere dense. We have $\cI \cdot \nu_*\cO_Y \subset \cO_X$. Let
$\cJ = \pi^* \cI \subset \cO_Z$. As $\cJ\cdot \cA^{\otimes\ell}$
is globally generated for some $\ell >0$, the linear system
$H^0(Z,\cJ\cdot \cA^{\otimes(\ell +1)}) \subset H^0(Z,
\cA^{\otimes(\ell +1)})$ embeds $Z\backslash V(\cJ)$ into some
projective space. Next, the multiplication with a canonical
section of $\cO_Z((\ell +1)\cdot E)$ defines a map $H^0(Z,
\cJ\cdot \cA^{\otimes(\ell +1)}) \to H^0(Z, \cJ\cdot
\pi^*\cL^{\otimes(\ell +1)}) \subset  H^0(Z,
\pi^*\cL^{\otimes(\ell +1)}) = H^0(X, \nu_*\nu^*
\cL^{\otimes(\ell+1)})$, whose composition with $H^0(X, \nu_*\nu^*
\cL^{\otimes(\ell+1)}) \to H^0(X, \cC\otimes
\cL^{\otimes(\ell+1)})$ is identically zero. So the image of
$H^0(Z, \cJ\cdot \cA^{\otimes(\ell +1)})$ in $H^0(Z, \cJ\cdot
\pi^*\cL^{\otimes(\ell +1)})$ is contained in the subspace
$H^0(X, \cL^{\otimes(\ell +1)})$. Hence global sections of
$\cL^{\otimes(\ell+1)}$ embed a Zariski open subset of $X$.
\end{proof}

We return to the situation of Theorem~\ref{criterion}. To show
that $L$ is big on the reduced space, we do not need that the
Lelong numbers of $h$ vanish. However, we give the more general
argument.

Let  $U \subset X_{reg}$ be a Zariski open subset, which is
quasi-projective. We can find a projective compactification $\ol
U$ together with a modification $\gamma: \ol U \to \ol X$ such
that the divisor $D = \ol U \backslash U$ has only normal
crossings singularities and such that the singular hermitian
metric $h$ extends from $U$ to $\ol U$ as a singular hermitian
metric on $\gamma^*L$ (cf.\ Section~\ref{singherm}. As usual one
can construct a complete K"ahler form $\eta_U$ on $U$ with
Poincar\'e growth near the boundary from a K"ahler form on $\ol U$
and a canonical section of $D$.

\begin{lemma}\label{nonvansec}
Let $x \in U$ be a point. Then there exists some $m_0>0$ so that for
any $m\geq m_0$ there is a section
$$
\sigma \in H^0_{(2)}(U, \cO_U(K_U+mL))
$$
with $\sigma(x)\neq 0$.
\end{lemma}
\begin{proof}
We use Kodaira's argument. Let $W=\{(z_1,\ldots, z_n)\}\subset U$
be a coordinate neighborhood and $\rho$ a cut-off function with
support in $W$, which is identically equal to one a relatively
compact neighborhood of $x$ contained in $W$ and has values
between $0$ and $1$. We set
$$
\psi_x= \rho(z) \cdot n\cdot\log(\sum|z_i|^2).
$$
There exists some $m_0>0$ and a continuous strictly positive
function $\wt c(p)$ on $U$ so that
$$
\sqrt{-1}\partial\ol\partial\psi_x + m_0\cdot \Theta_h \geq \wt
c(p) \cdot \eta_U.
$$
Let $m\geq m_0$. We chose a local section $\sigma_x \in
H^0(W,K_X+mL)$ on $W$ with $\sigma_x(x)\neq 0$ and set
$$
f=\ol\partial(\rho\sigma_x).
$$
By assumption that $h$ has absolutely continuous curvature on
$U$, so the Lelong numbers vanish on all curves, and by the
argument from Section~\ref{mulshea} the multiplier ideal sheaf is
equal to the structure sheaf, and $h$ is in $L^1_{loc}(U)$.
Moreover $f$ vanishes identically in a neighborhood of $x$. So
$$
\int_U \frac{1}{\wt c(p)} {\bf e}^{-\psi_x} ||f||^2 dV_\eta <
\infty,
$$
and Proposition~\ref{hoer} implies the existence of a
$(n,0)$-form $u$ of class $C^\infty$ (since $f$ is of class
$C^\infty$) with values in $mL$ such that
$$
\ol\partial u =f,
$$
and
$$
(\sqrt{-1})^{n^2} \int_U {\bf e}^{-\psi_x} h^m u \wedge \ol u <
\infty.
$$
The finiteness of the above integral (plus the fact that $h$ is
bounded from below, and that $u$ is holomorphic on some
neighborhood of $x$) imply that $u(x)=0$. Now we can see that
$$
\sigma = \rho \cdot \sigma_x -u
$$
is an element of $H^0_{(2)}(U, \cO_U(K_U + mL))$, which does not vanish
at $x$.
\end{proof}

In a similar way, by taking two points and directions at a point we
obtain the following lemma.

\begin{lemma}\label{em} For any compact set $K\subset U$ there exists
a number $m(K)>0$ so that for all $m \geq m(K)$ the linear system
$|H^0_{(2)}(U, \cO_U(K_U+mL))|$ gives an embedding of $K$ into a
projective space.
\end{lemma}

We have the following extension property:

\begin{lemma}\label{l2ext} There is a canonical embedding:
$$
H^0_{(2)}(U, \cO_U(K_U + mL)) \hookrightarrow H^0(\ol U, \cO_{\ol
U}(K_{\ol U}+ m \gamma^* L))
$$
\end{lemma}

\begin{proof} In our situation $h$ possesses an extension $\wt h$ as a
singular metric on $\gamma^* L$ with positive curvature in the
sense of Definition~\ref{singdef}. In particular, $\wt h$ is
locally bounded from blow by a positive constant.
 For any $\sigma \in
H^0_{(2)}(U, \cO_U(K_U + mL))$ we have
$$
(\sqrt{-1})^{n^2} \int_{\ol U} \wt h^m \gamma^*\sigma \wedge
\ol{\gamma^*\sigma} < \infty.
$$
So $\gamma^* \sigma$ extends holomorphically to $\ol U$.
\end{proof}

We are given a singular hermitian metric on $L$ over the singular
complex space $\ol X$, which amounts to a singular hermitian
metric $h$ on $\ol X_{reg}$, which can be extended from $U$  as a
singular metric $\wt h$ on $\gamma^* L$ over $\ol U$. The latter
defines a multiplier ideal sheaf $\cI(\wt h^m) \subset \cO_{\ol
U}$ which is defined by the following property: For all open
subsets $W \subset \ol U$ the space
$$
(\cO_{\ol U}(K_{\ol U} + m \gamma^* L)\otimes \cI(\wt h^m))(W)
$$
consists of all
$$
\sigma \in \cO_{\ol U}(K_{\ol U} + m \gamma^* L)(W \cap U)
$$
such that
$$
(\sqrt{-1})^{n^2} \int_V \wt h^m \sigma \wedge \ol \sigma < \infty
$$
for all $V \subset\subset W $.

\begin{definition}\label{bigherm}
A bundle $(L,h)$ is called {\em big} in the sense of singular
hermitian bundles, if
$$
\limsup_{m \to \infty} m^{-n} h^0(\ol U, \cO_{\ol
U}(m\gamma^*L)\otimes \cI(\wt h^m))
> 0
$$
holds.
\end{definition}
For any such bundle, the pull-back $\wh L$ of $L$ to the
normalization $\wh X \to \ol X$ satisfies $ \limsup_{m \to \infty}
m^{-n} h^0(\wh X,\cO_{\wh X}(m \wh L))>0$. According to
Proposition~\ref{equivbig}, any such bundle is big in the usual
sense, and Proposition~\ref{ZOemb} guarantees that associated
linear systems embed certain Zariski open subsets.

We claim:
\begin{proposition}\label{big}
The above line bundle $(L,h)$ is big.
\end{proposition}

\begin{proof} Let
$$
0 \neq \sigma_0 \in H^0_{(2)}(U, \cO_{\ol X}(K_{\ol X}+ m_0 L)),
$$
be a section, which we extend to $\ol U$. We denote by $D_0 \subset \ol
U$ the zero divisor. Next, we consider the restriction morphism
\begin{equation*}
\begin{split}
r_m: H^0_{(2)}(U, \cO_{\ol X}(K_{\ol X}+ m L)) =
H^0(\ol U, \cO_{\ol U}&(K_{\ol U} + m \gamma^* L ) \otimes \cI(\wt h^m))\\
\to  & H^0(D_0,\cO_{D_0}(K_{\ol U} + m \gamma^*L) ).
\end{split}
\end{equation*}
Since
$$
h^0(D_0, \cO_{D_0}(K_{\ol U} + m \gamma^* L )) = O(m^{n-1}),
$$ and
$$
\limsup_{m \to \infty} m^{-n} \dim H^0_{(2)}(U,\cO_{\ol X}(K_{\ol
X} + m L) ) > 0,
$$
we see that
$$
\limsup_{m \to \infty} m^{-n} \dim\ker  r_m > 0.
$$
Now
$$
 \ker  r_m \subset H^0(\ol U, \cO_{\ol U}(K_{\ol U} + m \gamma^* L)
\otimes \cI(\wt h^m)) \cap H^0(\ol U, \cO_{\ol U}(m-m_0) \gamma^*
L).
$$
Since $\wt h$ is locally bounded from below by some (positive)
constant (in the appropriate measure theoretic sense) $\cI(\wt
h^m)\subset \cI(\wt h^{m-m_0})$ holds. So
$
\ker  r_m \subset H^0(\ol U, \cO_{\ol U}((m-m_0) \gamma^* L)
\otimes \cI(\wt h^{m-m_0}))
$.
\end{proof}
We state the following general fact, which implies that the above
line bundle $L$, pulled back to a desingularization is nef.
\begin{proposition}\label{nef} Let $Y$ be a projective manifold
and $(L,h)$ a positive, singular hermitian line bundle, whose
Lelong numbers vanish everywhere. Then $L$ is nef.
\end{proposition}
\begin{proof}
Let $A$ be an ample line bundle on $X$. For any $y \in Y$ one
considers a finite, locally free resolution
$$
\cP^\bullet \to \mathfrak m_{Y,y}
$$
of the maximal ideal at $y$. Then we chose a multiple
$\ell(y)\cdot A$ so that all $\cP^j \otimes K_Y^{-1}(\ell(y)\cdot
A)$ are positive. The value for $\ell(y)$ can be taken uniformly
in a neighborhood of $y$. So we choose $\ell_0$ uniformly on $Y$
with this property. As the Lelong numbers vanish, the multiplier
ideal sheaves $\cI(h^m)$ are equal to $\cO_Y$ for all $m>0$. So
$H^q(Y,\cP^r \otimes \cO(\ell_0 A + m L))=0$ for all $r$ and
$q,m>0$ by the Nadel vanishing theorem. Now $H^1(Y,\mathfrak
m_{Y,y}\otimes \cO_Y(\ell_0 A + mL))=0$ for $m>0$, and the sheaves
$\cO_Y(\ell_0 A + mL)$ are globally generated. So $\cO_Y(\ell_0 A
+ m L)$ is globally generated, in particular nef. Hence $ L+
\frac{\ell_0}{m}A$ is nef for all $m>0$. With $m\to \infty$ the
claim follows.
\end{proof}

\section{Embedding of non-reduced spaces}\label{embnonred}
\begin{proposition}\label{embnonredsimp} Let $X$ be a compact
complex space, which possesses a holomorphic line bundle $L$,
whose restriction to the reduction $X_{red}$ is ample. Then $L$
is ample.
\end{proposition}

\begin{proof} Let $\cO_{X_{red}}= \cO_X/\cI$, and $X_{j}=(X_{red},
\cO_X/\cI^{j+1})$ so that $X_{red}=X_0$ and $X=X_k$ for some $k$.
Let $\cL=\cO_X(L)$, and  $\cL_0=\cL|X_0$. Now
$$
{\rm (E_j)}\hspace{1cm} 0 \to \cF_j \to \cO_{X_{j+1}} \to
\cO_{X_j} \to 0
$$
where $\cF_j=\cI^{j+1}/\cI^{j+2}$ is a coherent
$\cO_{X_0}$-module.

We can assume from the beginning that $\cL_0$ is very ample on
$X_0$, and that $H^1(X_0, \cF_j \otimes \cL^{\otimes\ell})=0$ for
all $\ell >0$, and $j=0,\ldots k-1$, furthermore that $\cF_j
\otimes \cL^{\otimes \ell}$ is globally generated for all $j$,
and all $\ell > 0$. Now for all $\ell>0$ the map $H^0(X,
\cL^{\otimes \ell}) \to H^0(X_0, \cL_0^{\otimes \ell})$ is
surjective. So we have a holomorphic map $\Phi_{|L|} :X \to
\mathbb P_N$, whose restriction $\Phi_0$ to $X_0$ is an
embedding, i.e $\cO_{\mathbb P_N} \to \cO_X$ followed by $\cO_X
\to \cO_{X_0}$ is surjective. We show by induction over $j$ the
existence of compatible maps $\alpha_j:\cO_{\mathbb P_{N(j)}} \to
\cO_{X_j}$. (In each step the number $N$ will have to be raised.)

We assume that $\cO_{\mathbb P_N} \to \cO_{X_j}$ is surjective.
The pull-back of the small extension $(E_j)$ of sheaves of
analytic $\mathbb C$-algebras with respect to this map induces the
direct sum of spaces $X_{j+1}\oplus_{X_j}\mathbb P_N$. As the
surjective morphism  $\alpha_j$ can be lifted to the morphism
$\alpha_{j+1}$, the induced small extension is trivial: we have
the following diagram.
$$
\begin{CD}
0 @>>> \cF_j @>>>\cO_{X_{j+1}} @>>> \cO_{X_j} @>>> 0 \\
@. @|@A{\rm surj}AA @A{\rm surj}A{\alpha_j}A @. \\
0 @>>> \cF_j @>>> \cO_{\mathbb P_N}[\cF_j] @>>> \cO_{\mathbb P_N}
@>>> 0
\end{CD}
$$
Denote by $i: X_{j+1} \hookrightarrow \mathbb P_N[\cF_j]$ and
$\tau:\mathbb P_N[\cF_j] \to \mathbb P_N$ resp.\ the embedding
and projection resp. Then $i^*\tau^*\cO_{\mathbb P_N}(1) = \cL$.
So in
$$
0 \to \cF_j \otimes \cL_0 \to \left( \cO_{\mathbb
P_N}[\cF_j]\right)\otimes\cO_{\mathbb P_N}(1) \to \cO_{\mathbb
P_N} \to 0
$$
we can identify the middle term with $\cO_{\mathbb P_N}(1)[\cF_j
\otimes \cL_0]$. Let $\{\sigma_0,\ldots,\sigma_N\} \subset
H^0(\mathbb P_N,\cO_{\mathbb P_N}(1))$ be a basis, and let
$\cF_j\otimes \cL_0$ be generated by global sections
$\tau_q,\ldots,\tau_r$. Let $\epsilon^2=0$. Then the
$\sigma_0,\ldots, \sigma_N, \epsilon \tau_1, \ldots,
\epsilon\tau_r$ give rise to an embedding $\mathbb P_N[\cF_j]
\hookrightarrow \mathbb P_{N+r}$. Altogether $\cO_{\mathbb
P_{N+r}} \to \cO_{X_{j+1}}$ is surjective. \end{proof}

We need the above statement in a more general situation.

Let $\ol Z$ be a non-reduced complex space equipped with a
holomorphic line bundle $L$, $\cL=\cO_{\ol Z}(L)$. Let $\ol X =
{\rm red}(\ol Z)$, and let $X \subset \ol X$ be a Zariski open
subset. We denote by $Z$ the restriction of the non-reduced
structure to $X$. We consider the meromorphic map
$\Phi=\Phi_{|L|_{\ol X}|}: \ol X \hookrightarrow \mathbb P_N$.
\begin{proposition}\label{nonreduced}
Assume that $\Phi|_X :X \to \mathbb P_N$ is an embedding, and let
$L$ be also nef. Then for some multiple $\ell_0$ the meromorphic
map $\Phi_{|\ell_0L|}: \ol Z \to \mathbb P_M$ defines an
embedding of an open subspace of $Z$.
\end{proposition}
\begin{proof}
First, we take a (projective) desingularization of $\ol X$, and
pull back the meromorphic map $\Phi$. Then we eliminate the
indeterminacy set by a sequence of blow-ups with smooth centers.
This procedure is locally done by embedding the space in a smooth
ambient space, blowing up the ambient space along smooth centers
several times, and by taking in each step the proper transform of
$\ol X$. We take locally embeddings of $\ol X$, which extend to
embeddings of $\ol Z$. However, in general there exist no proper
transforms of non-reduced spaces. So in all steps, we take the
full pre-image under the monoidal transformation, and arrive at a
space $\rho^\vee: \wt Z^\vee \to \ol Z$ together with the
restriction $\pi : \wt X \to \ol X$, which allows a morphism
$\Psi : \wt X \to \mathbb P_N$. We denote by $\ol Z_k \subset \ol
Z$ the $k$-th infinitesimal neighborhood of the compact manifold
$\ol X$ in $\ol Z$ so that $\ol Z_0=\ol X$ and $\ol Z_{k_0} = \ol
Z$ for some $k_0>0$. We consider the $k$-th infinitesimal
neighborhoods $\wt Z_k$ of $\wt X$ in $\wt Z^\vee$, and we define
$\wt Z:=\wt Z_{k_0}$, and we denote the map to $\ol Z$ by $\rho$.
The infinitesimal neighborhoods give rise to short exact sequences
$$
0 \to \cF_j \to \cO_{\wt Z_{j+1}} \to \cO_{\wt Z_{j}} \to 0,
$$
where the $\cF_j$ are coherent $\cO_{\wt X}$- modules.

Let $n=\dim X$. Denote by $\wt L$ the pull-back of $L$ to $\wt Z$.
We claim that (for any fixed $j$) $h^1(\wt X, \cF_j(\ell \cdot \wt
L|_{\wt X})) O(\ell^{n-1})$: The bundle $L|_{\wt X}$ is big. After
replacing $\wt L$ by a multiple we write $\wt L|_{\wt X}= A + E$,
where $A$ is ample and $E$ is effective by Kodaira's lemma. As
$\wt X$ is projective and $L$ is nef and big, we can fix a K"ahler
form $\eta_{\wt X}$ and find hermitian metrics $h_r$ on $\wt L|\wt
X$ such that the curvature of $h_{r}$ is greater or equal to
$-(1/r)\eta_{\wt X}$ for all $r>0$. This shows the existence of
some $m_0$ such that
$$
H^1(\wt X, \cF_j \otimes \cO_{\wt X }(m A + \ell \wt L))= 0
$$
for all $m \geq m_0$ and $\ell>0$.

Let $mE$ denote the non-reduced space with support $E$, induced
by the divisor $mE$. Then
$$
0 \to \cF_j(m A + \ell \wt L) \to \cF_j((m+\ell)\wt L) \to
\cF_j((m+ \ell)\wt L)|_{mE} \to 0
$$
is exact, and for $m\geq m_0$
$$
0 \to H^1(\wt X, \cF_j((m+\ell)\wt L)) \to H^1(mE, \cF_j((m+\ell)
\wt L)|_{mE})
$$
as well.

We fix $m=m_0$ and look at $\ell \gg 0$. Then $h^1(\wt X,
\cF_j((m_0+\ell)\wt L))= O(\ell^{n-1})$. Now
\begin{gather*}
H^0(\wt X,\cO_{\wt X_{j+1}}((m_0+\ell)\wt L)) \to H^0(\wt
X,\cO_{\wt X_{j}}((m_0+\ell)\wt L))\\ \hspace{6cm} \to H^1(\wt X,
\cF_j((m_0+\ell)\wt L))
\end{gather*}
is exact, and $h^0(\wt X,\cO_{\wt X_{j}}(m_0+\ell)\wt L)$ grows
like $\ell^n$, because we can assume by induction that high
powers of $\wt L$ embed a Zariski open subset of $\wt X_j$, so
$h^0(\wt X,\cO_{\wt X_{j+1}}((m_0+\ell)\wt L)) \sim \ell^n$. This
means that the sections of $H^0(\wt X,\cO_{\wt
X_{j+1}}((m_0+\ell)\wt L))$ define a meromorphic map $\Xi: \wt
X_{j+1} \to \mathbb P_M$, which embeds an open subset $W\cap \wt
X$, $W\subset \wt X_{j+1}$ of the reduction $\wt X$(cf.\
Proposition~\ref{ZOemb}). We assume that $W$ is closed in some
open set $\mathbb P_M \backslash A$. The sets $A$ and $W$ can
also be chosen in a way that
$$
0 \to \Xi_* \cF_j|\mathbb P_M\backslash A  \to \Xi_*\cO_{\wt
Z_{j+1}}|\mathbb P_M\backslash A \to \Xi_*\cO_{\wt Z_{j}}|\mathbb
P_M\backslash A \to 0
$$
is still exact. As in the proof of
Proposition~\ref{embnonredsimp} we consider the fibered sum of
complex spaces $W\oplus_{W\cap \wt X_j} (\mathbb P_M\backslash
A)$, which is isomorphic to the trivial extension $(\mathbb
P_M\backslash A)[\cF_j|W]$, which is clearly quasi-projective.
The rest follows like in the proof of Proposition~\ref{ZOemb}.
\end{proof}

\section{Proof of the quasi-projectivity criterion}
We first need the following fact:
\begin{lemma}\label{updown}
Let $\pi: Y\to X$ be a proper holomorphic map of reduced, compact
not necessarily normal complex spaces. Let $S\subset X$ be a
closed subspace such that $\pi$ is an isomorphism over
$X\backslash S$. Let $\cI=\cI_S \subset \cO_X$ be the vanishing
ideal of $S$. Then for any coherent sheaf $\cF$ on $X$ there is a
number $m >0$ and morphism $\mu: \cI_S^m \cdot (\pi_*\pi^*\cF) \to
\cF$, which is an isomorphism over $X\backslash S$.
\end{lemma}
\begin{proof}
We consider the short exact sequence $0\to \cF \to \pi_*\pi^* \cF
\to \cC \to 0$, where ${\rm supp}(\cC)\subset S$. Now the zero
set of the annihilator ideal $V(Ann_{\cO_X}(\cC))$ is contained
in $S$ so that $\cI^m \cdot \cC=0$ for some $m>0$.
\end{proof}

We consider the situation of Theorem~\ref{criterion} and assume
that $\ol X$ is reduced and irreducible. Let
$$
S= \{x \in X ; |mL| {\rm \;does\; not\; define\; an \;
 embedding \; around\;}  x {\rm \; for \; all\; } m > 0
\}.
$$

From Proposition~\ref{big}, we know that the line bundle $L$ in
$\ol X$ is big, and by Proposition~\ref{ZOemb} the linear system
$|mL|$ provides an embedding of some Zariski open subspace. Using
Noether induction, we see that there is a number $m_0>0$ such
that $\Phi_{|m_0L|}$ embeds $X\backslash S$. In particular, it
embeds $X$, if $S$ is empty.

\begin{lemma}\label{glgen} Let $\cF$ be a coherent $\cO_{\ol
X}$-module. Then there exists some $\ell_0 > 0$ such that for all
$\ell \geq \ell_0$ the sheaf $\cF \otimes \cO_{\ol X}( \ell m_0
L)$ is generated by global sections at all points $x \in X
\backslash S$.
\end{lemma}

\begin{proof}
Let $\Phi=\Phi_{|m_0L|}$, and denote the graph of $\Phi$ by
$\Gamma_\Phi$. We have a diagram
\begin{equation*}
\settriparms[1`1`1;350] \btriangle[\Gamma_\Phi`\ol X` \mathbb
P_N;\pi`\Psi`\Phi].
\end{equation*}
Let $\wt S = S \cup (\ol X \backslash X)$, $T= \Psi(\pi^{-1}\wt
S)$, and $\cI_T \subset \cO_{\mathbb P_N}$ the corresponding
ideal. According to Serre's theorem, for any $m_1$ the sheaf $
\cI_T^{m1} \cdot \Psi_* \pi^* \cF \otimes \cO_{\mathbb
P_N}(\ell)$ is generated by global sections for all $\ell
\allowbreak \geq \allowbreak \ell_0(m_1)\allowbreak >\allowbreak
0$.

From the construction, we have a morphism of sheaves $\Psi^*
\cO_{\mathbb P_N}(1) \to \pi^* \cO_{\ol X}(m_0 L)$, which is an
isomorphism over $\Gamma_\Phi\backslash \pi^{-1}(\wt S)$. We have
the following morphisms.
\begin{gather*}
 \cI_T^{m_1}\cdot
(\Psi_*\pi^*\cF) \otimes \cO_{\mathbb P_N}(\ell) \to
\cI_T^{m_1}\cdot \Psi_*\pi^*(\cF \otimes \cO_{\ol X}(m_0\ell L ))
\\ \hspace{5cm} \to \Psi_*(\cI_{\pi^{-1}\wt S}^{m_1}\cdot \pi^*(\cF \otimes \cO_{\ol
X}(m_0\ell L )))
\end{gather*}
Now
\begin{gather*}
H^0(\mathbb P_N, \cI^{m_1}\cdot(\Psi_*\pi^*\cF) \otimes
\cO_{\mathbb P_N}(\ell )) \to H^0(\Gamma_\Phi, \cI_{\pi^{-1}\wt
S}^{m_1}\cdot \pi^*(\cF \otimes \cO_{\ol X}(m_0\ell L )))\hspace{1cm}\\
\hspace{1cm} = H^0(\ol X, \pi_*(\cI_{\pi^{-1}\wt
S}^{m_1}\pi^*\cF)\otimes \cO_{\ol X}(m_0\ell L)) \to \hspace{6cm}\\
H^0(\ol X, \cI_{\wt S}^{m_2}\cdot \pi_*\pi^* \cF \otimes \cO_{\ol
X}(m_0 \ell L)) \to H^0(\ol X, \cF \otimes \cO_{\ol X}(m_0 \ell
L)).
\end{gather*}
Here, we chose $m_1\gg 0$ large enough so that
a corresponding $m_2$ satisfies $\cI_{\wt S}^{m_2}\pi_*\pi^*\cF
\subset \cF$.

Over $X\backslash S$ the above morphism of sheaves are
isomorphisms so that we can produce enough global sections, which
generate  $\cF \otimes \cO_{\ol X}(m_0\ell L)$ over $X\backslash
S$.
\end{proof}

In the above situation we also need the case, where $\cF$ is an
ideal in $\cI \subset\cO_{\ol X}$. If $\ol X$ is smooth or
normal, we have automatically $\pi_*\pi^* \cI = \cI$.

\begin{definition}
Let $Y$ be a reduced complex space of pure dimension $n$. The
$L^{2}$-dualizing sheaf $\omega^{(2)}_{Y}$ of $Y$ is defined by
\begin{equation*}
\begin{split}
\omega^{(2)}_{Y}&(W) = \{ \eta\in \Gamma
(W_{reg},{\cO}(K_{Y_{reg}})) ; \\
&\qquad \qquad(\sqrt{-1})^{n^2} \int_{V}\eta\wedge\bar{\eta} <
\infty \,\,\mbox{for every $V \subset\subset W$}\},
\end{split}
\end{equation*}
where $W$ runs through the open sets of $Y$.
\end{definition}

If $\alpha: \wt Y \to Y$ is a desingularization such that the
singular locus of $Y$ corresponds to a normal crossings divisor
in $\wt Y$, we have $\omega^{(2)}_{Y}=\alpha_*\cO_Y(K_Y)$. In
particular, $\omega^{(2)}_{Y}$ is coherent.

Now we set
$$
S_+ = \{x \in X ; \omega^{(2)}_{\ol X} \otimes L^{\otimes m}
\text{\it\/ is not generated by global sections at\, }x
\text{\it\, for all\, } m \},
$$
and
$$
S_- = \{x \in X ; \omega^{(2)\vee}_{\ol X} \otimes L^{\otimes m}
\text{\it\/ is not generated by global sections at\, }x
\text{\it\, for all\, } m \}.
$$
Now, by definition $S_+\cup S_- \supset S$, and the converse
inclusion follows from Lemma~\ref{glgen}, so
$$
S_+\cup S_- = S.
$$
Let $R$ be the non-normal locus of $X$. We denote by $\cI_{\wt S
\cup R} \subset \cO_{\ol X}$ the ideal of functions that vanish
on $\wt S \cup R$. Lemma~\ref{glgen} implies that there exists
$m_1
> 0$ such that $\cO_{\ol X}(m_0m_1 L)\otimes \cI_{\wt S\cup R}$ is
generated by global sections at all points $x \in X \backslash
S$. Let $\{\sigma_0, \ldots, \sigma_{N(m_0m_1)} \}$ be a $\mathbb
C$-basis of $\Gamma(\ol X, \cO_{\ol X}(m_0m_1 L)\otimes \cI_{\wt
S\cup R})$. Then
$$
h_0= \frac{1}{\sum |\sigma_i|^2}
$$
defines a singular hermitian metric on $m_0m_1L$ over the space
$\ol X$, whose singularities are contained in $\wt S\cup R$. Let
$m_+$ and $m_-$ resp.\ be integers such that $\omega^{(2)}_{\ol X}
\otimes \cO_{\ol X}(m_+m_0 L)$ and $\omega^{(2)\vee}_{\ol X}
\otimes \cO_{\ol X}(m_- m_0L)$ are generated by global sections
over $X\backslash S$. Let $\{\sigma_{+,k}\}$ be a basis for the
space of sections of the former sheaf over $\mathbb C$. We define
a hermitian metric $h_+$ on $\omega^{(2)}_{\ol X} \otimes \cO_{\ol
X}(m_+ m_0L)$ by
$$
h_+=  \frac{1}{\sum|\sigma_{+,k}|^2}.
$$
In a similar way a metric $h_-$ on $\omega^{(2)\vee}_{\ol X}
\otimes \cO_{\ol X}(m_-m_0 L)$ is constructed. These metrics are
so far only well-defined, where the sheaves are locally free,
i.e.\ at least on $(\ol X)_{reg}$. We chose a desingularization
$\pi : \wt X \to \ol X$ in such a way that also
$\pi^*\omega^{(2)}_{\ol X}$ is invertible (after dividing by the
torsion part). Then the pull-backs of the sections $\sigma_{+,k}$
and $\sigma_{-,k}$ define singular hermitian metrics over $\wt X$
on the corresponding line bundles. We impose a further condition:
Let $U \subset X_{reg}$ be a Zariski open subset, which is
quasi-projective. Let $\ol U$ be a projective compactification
that dominates $\wt X$ with modifications $\rho : \ol U \to \wt X$
and $\gamma : \ol U \to \ol X$. We pull back $h$, $h_+$, and $h_-$
back to $\ol U$, and we assume as above that $\ol U \backslash U $
is a divisor with normal crossings singularities. We denote by
$\eta_{\ol U}$ a K"ahler form on $\ol U$ and by $\eta_U$ a complete
K"ahler form on $U$ with Poincar\'e growth condition near the
boundary as above. By Proposition~\ref{big} the line bundle
$\gamma^*L$ is big on $\ol U$. Kodaira's lemma provides an
effective $\mathbb Q$-divisor $A$ such that the $\mathbb
Q$-divisor $\gamma^* L - A$ is ample, giving rise to a strictly
positive hermitian metric $h'$ of class $C^\infty$ on the $\mathbb
Q$-line bundle $\gamma^*L -A$. Let $a\cdot A$ be a divisor with $a
\in \mathbb N$, and $\sigma_A$ a section of $\cO_{\ol U}(a\cdot
A)$. Then
$$
\frac{h'}{|\sigma_A|^{2/a}}
$$
defines a singular hermitian metric on $\gamma^* L$, whose
curvature $\Theta$ satisfies
$$
\Theta \geq \alpha \cdot \eta_{\ol U}
$$ for some $\alpha >0$ on $\ol U$. Let $D_j$ be the components of
the normal crossings divisor $D=\ol U \backslash U$. We equip the
bundles $[D_j]$ with a $C^\infty$ hermitian metric. We can find
canonical sections $\tau_j$ and some $\beta > 0 $ such that the
curvature of the modified hermitian metric
$$
h''= \frac{h'}{|\sigma_A|^{2/a}} \cdot \prod \left( -\log
\|\tau_j\|\right)^\beta
$$
over $U$ satisfies
$$
\Theta_{h''} \geq \varepsilon \cdot \eta_{\ol U}.
$$

The following considerations apply to the above line bundles on $\ol
U$.

Let $p,r \in \mathbb N$, and let $1>\delta > 0$. Then
$$
\wh h:= h_{p,r,\delta}:= \gamma^*( h_0^{p-\delta}h_+^{r}
h_-^{r+1}) {h''}^{\; \delta m_0m_1}
$$
is a singular hermitian metric on
$$
\gamma^*(L^{\otimes e}\otimes \omega^{(2)\vee}_{\ol X} )
$$
with $e=m_0(pm_1 + rm_+ + (r+1)m_-)$ and $\Theta_{\wh h}\geq
\varepsilon \delta m_0m_1 \cdot \eta_{\ol U}$. Hereafter we shall
consider $\wh h$ as a singular hermitian metric on $L^{\otimes
e}\otimes \omega_{\ol X}^{(2)\vee}$.

We choose $p>0$ large enough so that $\cI(h^p_0)$ annihilates
$\pi_*\cO_{\ol U}/\cO_{\ol X}$.

Although as a singular hermitian metric on a {\em line bundle}
$\wh h $ is only defined over ${\ol X}_{reg}$, a coherent
multiplier ideal sheaf $\cI(\wh h )\subset \cO_{\ol X}$ can be
given a meaning as follows: For $W \subset \ol X$ open, we define
\begin{equation*}
\begin{split}
(\cI(\wh h ) \otimes \cO_{\ol X}(L^{\otimes e}))(W)&=\\
\{\sigma \in \Gamma(W,\cO_{\ol X}&(L^{\otimes e})) ; \int_{V\cap
(\ol X)_{reg}} |\sigma|^2 \wh h  < \infty
\text{\it \; for all } V \subset\subset W \}.\\
\end{split}
\end{equation*}
Observe that $\wh h |(\ol X)_{reg}$ is a (singular) hermitian
metric on $(L^{\otimes e}|_{(\ol X)_{reg}})\otimes K^{-1}_{(\ol
X)_{reg}}$. The ideal is $\cI(\wh h )$ is coherent, since
\begin{equation}\label{ideal}
\cI(\wh h )\otimes \cO_{\ol X}(L^{\otimes e}) = \cI(\gamma^*\wh
h)\cdot \gamma_*(K_{\ol U}\otimes \gamma^*(L^{\otimes e} \otimes
\omega^{(2)\vee}_{\ol X}) )
\end{equation}
holds by the usual definition  of the usual multiplier ideal sheaf
$\cI(\gamma^*\wh h)$ for the singular hermitian metric
$\gamma^*\wh h$ on $\gamma^*(L^{\otimes e}\otimes
\omega^{(2)\vee}_{\ol X})$.

Now we determine the value of $\delta > 0$.

\begin{remark}
We can see from the definition of $\wh h $ on $\ol X$ that for
sufficiently large $p$ the zero set  $V(\cI(\wh h ))\cap X = S$.
Furthermore for large $p$ the embedding dimension of the
non-reduced space defined by $\cI(\wh h )$ is equal to $\dim \ol
X$.
\end{remark}
\begin{proof}  For large $p$ the
contribution of $h_0^{p-1}$ to $\gamma^*( h_0^{p-\delta}h_+^r
h_-^{r+1})$ dominates the rest, so that the zero set of the
multiplier ideal sheaf is contained in $\wt S\cup R$ and contains
$\ol S \cup R$. Next $\delta >0$ is chosen small enough: The term
$1/|\sigma_A|^2$ is equipped with a small exponent so that the
$L^2$-integrability condition for holomorphic sections is not
affected, and $V(\cI(\wh h))\cap X = S \cup R$ still holds. For
large $p$ also the second statement is satisfied.
\end{proof}

\begin{proposition}\label{extfrominf} The canonical map
$$
H^0(\ol X, \cO_{\ol X}(L^{\otimes e})) \to H^0\left(\ol X,
\cO_{\ol X}(L^{\otimes e}) \otimes \big(\cO_{\ol X}/\cI(\wh h
)\big)\right)
$$
is surjective.
\end{proposition}

\begin{proof} Let $\tau \in H^0(\ol X, \cO_{\ol X}(L^{\otimes e})
\otimes (\cO_{\ol X}/\cI(\wh h)))$ be a section. For any
neighborhood $\wt W$ of $\wh S= V(\cI(\wh h))$ we can find a
$C^\infty$ section $\wt \tau$ of $L^{\otimes e}$, whose
restriction to $(\wh S, \cO_{\ol X}/\cI(\wh h)) $ equals $\tau$
with $\text{supp}(\wt \tau) \subset \wt W$. We consider
$\gamma^*\ol
\partial \wt \tau = \ol \partial \gamma^* \wt \tau$ on $\ol U$.

Since $H_{(2)}^1(U, \gamma^*(L^{\otimes e}\otimes
\omega^{(2)\vee}_{\ol X}) \otimes \cO_{\ol U}(K_{\ol U})) =0$ by
Nadel's vanishing  theorem (or H"ormander's theorem on
$L^2$-estimates resp.), there exists a $C^\infty$-section $u$ of
$\gamma^*(L^{\otimes e}\otimes \omega^{(2)\vee}_{\ol X}) \otimes
\cO_{\ol U}(K_{\ol U})$ on $U$, which is square-integrable with
respect to the singular hermitian metric and the complete K"ahler
metric on $U$ such that $\ol \partial u = \ol \partial \gamma^*
\wt \tau$.

So $0= \ol \partial(\gamma^* \wt\tau - u)$, i.e.\ $v = \gamma^*
\wt\tau - u \in H^0(U, \gamma^*(L^{\otimes e}\otimes
\omega^{(2)\vee}_{\ol X}) \otimes \cO_{\ol U}(K_{\ol U}))$.

We claim that $v$ is square-integrable: By (\ref{ideal}) $\| u
\|^2 \wh h$ is integrable over $(\ol X)_{reg}$.  Since $\wh h$ is
a singular hermitian metric of positive curvature, it is locally
bounded from below by a positive constant. Moreover $\wt \tau$ is
of class $C^\infty$, and $U$ carries the complete metric $\eta_U$
(with Poincar\'e growth condition). So $v$ extends
holomorphically to $\ol U$. Then $v$ gives rise to a holomorphic
section of $L^{\otimes e}$ on $\ol X$, which coincides with
$\tau$, when restricted to the subspace $(\ol S, \cO_{\ol X}/
\cI(\wh h ))$ (cf.\ equation (\ref{ideal})).
\end{proof}

{\it Proof of the Theorem.} We know from the induction hypothesis
that $m_0 L|\ol S$ defines an embedding of $S$. In particular,
$m_0 L|\ol S$ is big. Now $e$ is a sufficiently high multiple of
$m_0$, and also $e L|\ol S$ embeds $S$. In the last step, we need
to raise the power of $L$ without affecting the multiplier ideal
sheaf $\cI(\wh h)$. Here, we replace $L^{\otimes e}$ by
$L^{\otimes (e+t)}$, where the factor $L^{\otimes t}$ is equipped
with the original metric $h^t$. Since the curvature of $h$ is
absolutely continuous, by Proposition~\ref{id_mult}, we may
assume that $\cI(\wh h) = \cI(\wh h\cdot h^t)$ holds over $X$, if
we perturb $\delta$ by a small amount (i.e.\ we perturb $\wh h$).
However, this metric is not of analytic singularity, but a
singularity of type $\log(-\log(t))$ is negligible (cf.\
Remark~\ref{perturb}). We chose $t$ large enough so that
$L^{\otimes(e+t)}$ defines an embedding of a Zariski open
subspace of $(\ol S, \cO_{\ol X}/ \cI(\wh h ))$ by
Proposition~\ref{nonreduced}. Now $L^{\otimes(e+t)}$ embeds
$X\backslash S$ as well as a non-empty open subset of $S$, and it
also separates normal directions of this set in $X$. This
contradicts the choice of $S$. This proves
Theorem~\ref{criterion} for reduced spaces $\ol X$.

For non-reduced spaces again we use induction over the dimension.
Let $L$ be a line bundle on $\ol X$ with  the above assumptions.
We know that for some $m>0$ the meromorphic map
$\Phi_{|mL|red(\ol X)|}$ embeds $red(X)$. By
Proposition~\ref{nonreduced} we can choose $m>0$ so that
$\Phi_{|mL|}$ embeds a Zariski open subspace $X'\subset X$. Let
$T=red(X)\backslash red(X')$. According to the above proof, a
multiple of $L$, restricted to a high infinitesimal neighborhood
$T_{inf}$ of $T$ gives rise to a linear system, which embeds a
Zariski open subspace of $T_{inf}$. Finally, by
Proposition~\ref{extfrominf} global sections over $T_{inf}$ can
be extended to all of $\ol X$ contradicting the choice of $T$.
\qed

Finally Theorem~\ref{main1} and Theorem~\ref{criterion} imply
Theorem~\ref{main}.

\end{document}